\documentstyle[10pt]{amsart} %[8pt]
\begin{document}

\def\RR{\mathbb R}
\def\CC{\mathbb C}
\def\lb{\lambda}
\def\2{{1\over 2}}
\def\w{|w>}
\def\wl{|w>_{\lambda}}
\def\wwl{{}_{\lambda}<w|}
\def\ww{<w|}
\def\f{1\over}
\def\d{\partial\over\partial}
\def\vL{V_{\lambda}}
\def\vLL{V_{\lambda+1}}
\def\nn{\nonumber}
\newcommand{\De}{\Delta}

\def\lD{{\Large Definition }}
\def\lL{{\Large Lemma }}
\def\lT{{\Large Theorem }}
\def\lPn{{\Large Proposition }}
\def\lC{{\Large Corollary  }}
\def\lP{{\Large Proof: }}
\def\QD{{\Large Q.E.D.}}
\def\Rem{{\it{ Remark: }}}

\def\bea{\begin{eqnarray}}
\def\eea{\end{eqnarray}}
\def\bqa{\begin{eqnarray}}
\def\eqa{\end{eqnarray}}
\def\pr{\partial}
\def\apr{\overline{\partial}}
\def\noi{\noindent}
\def\nn{\nonumber}
\def\baselinestretch{1.5}
\def\beq{\begin{equation}}
\def\eeq{\end{equation}}
\def\ba{\beq\new\begin{array}{c}}
\def\ea{\end{array}\eeq}
\def\be{\ba}
\def\ee{\ea}

\begin{center}

\hfill ITEP-TH-16/99\\

\hfill \hfill {\tt math.RT/?????}\\

\bigskip
\bigskip
{\Large\bf  Raising operators for the Whittaker wave functions of the Toda
chain and intertwining operators.\\}

\bigskip
\bigskip
{\large A.Chervov \footnote{E-mail: chervov at vitep1.itep.ru, ~~~~
alex at lchervova.home.bio.msu.ru
} } \\
%\\} } \\
\medskip
{\it Moscow State University\\}
\medskip{\it and \\}
\medskip
{\it Institute of Theoretical and Experimental Physics\\
117259 Moscow, Russia\\}
\bigskip
\bigskip
\bigskip
%{Revised \today}
{\large Abstract.}
\end{center}

 Intertwiners between  representations of Lie groups can be used to
obtain relations for matrix elements. We apply this technique
to obtain different identities for the wave functions of the
open  Toda chain, in particular
raising operators  and bilinear
relations for the wave functions at different energy levels.
We also recall the  group theory approach to the Toda chain:
treating the wave functions as
matrix elements in irreducible
representations between the so-called Whittaker vectors, integral
representations of the wave functions, etc.

%%%%%%%%%%%%%%%%%%%%%%%%%%%%%%%%%%%%%%%%%%%%%%%%%%%%%%%%%%%%%%%%%%%
%%%%%%%%%%%%%%%%%%%%%%%%%%%%%%%%%%%%%%%%%%%%%%%%%%%%%%%%%%%%%%%%%%%
\section{Introduction}\label{intro}
%%%%%%%%%%%%%%%%%%%%%%%%%%%%%%%%%%%%%%%%%%%%%%%%%%%%%%%%%%%%%%%%%%%
%%%%%%%%%%%%%%%%%%%%%%%%%%%%%%%%%%%%%%%%%%%%%%%%%%%%%%%%%%%%%%%%%%%

  Interrelations between the theory of integrable systems and representation
group theory proved to be fruitful for the both fields. The subject has become
so vast and diverse that it looks  impossible to pick up a few of surveys
which would cover all the important questions. Some topics are covered in reviews
\cite{Kyoto,S,Dr-Sok} and books \cite{B,EFK}. Review \cite{Mir} is  closely
related to our consideration. Usually it is possible to find a group
theory interpretation for all the important structures of
integrable systems. Phase space, symplectic structure, Hamiltonians, Lax
pairs, etc.  can be interpreted in terms of an underlying group symmetry for
both  classical and quantum integrable systems.  However we will not
discuss these questions here.

%   Integrable systems, since the pionering work \cite{Krus} being the
%% starting point of the modern era of the research in this field,
%% became complicated
%% and diverse field of the modern mathematics and physics.
%% Recent appearance of the soliton equations in such different fields as
%% condensed matter: quantum Hall effect; quantum field theory  :
%% Seiberg-Witten theory \cite{MMM1}; algebraic geometry: quantum cohomology
%% \cite{Dubr}, \cite{Givent},\cite{Eguchi}, Langlands program \cite{Frenkel}
%% shows that the story of intergability is far from being finished.
%% Phenomena of integrability is closely related to the phenomena of
%%%the hidden symmetry, which is the real reason for the system to be integrable.
%%%Exploring the intergable systems from the point of view of the symmetry group
%%%has its own long story, one can look at \cite{OP} for the survey of the
%%%activity at seventies at this field, but it seems that now
%%%one cannot name 3-4 surveys which covers all important aspects of
%% the interrelations between group theory and integrable systems.
%% Some aspects are covered in surveys: \cite{Drinf-Sokol},
%% \cite{Kyoto}, \cite{JM}, \cite{Frenkel}, \cite{Gorsky}, \cite{Kharchev}, \cite{Mironov},
%% \cite{Morozov1}, \cite{Morozov2},
%% %\cite{Kirillov},
%% \cite{Cherednik},
%% \cite{Verdier1, Verdier2} and books \cite{Presley-Segal}, \cite{Kac},
%% \cite{Frenkel-Eti-Kir},\cite{Lezn-Sav}
%%%%\cite{Varch}.
%%%

  The main goal of this paper is to demonstrate how representation theory,
and, more concretely,
 {\it intertwining operators, can be used to obtain relations
between wave-functions at different energy levels.}
We will work out only the simplest case of the open Toda chain with two and three
particles, but the method can be applied
to other integrable systems.
The paper mainly follows the ideas of \cite{Liouvil,gklmm}.

  Let us recall one of the main ideas of the  group theory approach
to integrable systems: {\it solutions of integrable equations
usually are  matrix elements or traces of the groups elements in irreducible representations}.
Main examples we keep in mind are: the tau-function of the Kadomtsev-Petviashvili
which is a matrix element in fundamental representations of the group
$GL_{\infty}$ (see \cite{Kyoto,Mir,gklmm}), and the wave functions of
many-body problems like Toda or Calogero-Sutherland-Moser, which are matrix
elements in infinite-dimensional irreducible representations of the
semisimple groups (see for example \cite{Liouvil}).

%The idea particular examples of which was first
% presented in \cite{gklmm},\cite{Feig-Anton} and which we will elaborate here
The idea  to obtain different relations for wave functions is quite simple
and goes as
follows:
let functions $f_i$ $(g_i)$ be  solutions of some integrable system corresponding
to some representation $V_i$ ($W_i$ respectively ) of some group,
(it can be the trace or matrix element, can be a wave function or a tau function,
etc.)
{\it If one has an intertwinning operator
$V_1\otimes V_2\otimes ... \otimes V_n \to
W_1\otimes W_2\otimes ... \otimes W_m  $
then  there  'usually'{\footnote{'usually' means
that there are no such
relations for generic matrix elements,
while they typically emerge for solutions of
integrable systems.}}
 exist the following relation:
$ D (f_1 f_2 ... f_n) = \tilde D (g_1 g_2 ... g_m) $ },
where $D, \tilde D$ are some differential operators.
Miraculously enough this  simple idea leads to non-trivial relations.
It was demonstrated in \cite{gklmm,L-M} that in such way one
can obtain famous bilinear identities for KP,KdV tau-functions,
which corresponds to intertwiners between fundamental representations
of the $GL_{\infty}$, $\hat {SL}(n)$. % It was also proposed there
%how to write the tau-functions corresponding to nonfundamental
%represenation and to quantum deformation $U_q(sl_n)$, see \cite{Chervov}.
In \cite{Feig-Anton}, it was shown, that, in this way, one can obtain
the Baxter relation for transfer matrices of integrable lattice models.
The above idea is also used in representation theory approach
to Knizhnik-Zamolodchikov equation \cite{FR,TK,EFK},
and  in the theory of special functions
(see, for example \cite{Vilen,Jeugt}).

 %In this paper we will show how
 %to obtain this way different relations for Toda wave-functions.

  %Let us turn to the case considered
  In this paper,
we show how
to obtain different relations for the Toda wave-functions applying the above idea.
The group theory approach to the classical Toda chain has been first proposed
in \cite{Bogo}, group theory approach to the quantum Toda chain
has been proposed by Kostant and Kazhdan at late seventies (unpublished) and
further
elaborated in \cite{Semen}, \cite{Liouvil}. Let us recall it following
\cite{Liouvil}.  Consider the group $SL(n+1,\RR)$ with   the
standard Chevalley generators $e_i, h_i, f_i$, and let $V_\lambda$ be
irreducible representations with weight $\lb$ of the principal series
with weight $\lb$.
Let us introduce
the so called Whittaker vector
$|w> \in V_\lb$, i.e. the eigenvector of generators $e_i$.  Then,
one considers  the Whittaker function, i.e. matrix element $ W_\lb(\phi_1,
..., \phi_n)= <w| exp(\phi_1 h_1+ \phi_2 h_2+...+\phi_n h_n) |w>$.
The main observation is that {\it function $ W_\lb(\phi_1, ..., \phi_n)$
%called
%Whittaker function as function of the variables $ \phi_1, ..., \phi_n $
is a common eigenfunction of Toda Hamiltonians, which turn
out to be Casimirs of the  group $SL(n+1,\RR)$}, while the weight $\lb$
marks the corresponding energy level.
%The Whittaker vector is the following: it is such vector
%in $V_\lb$ such that it is eigenvector for all generators $e_i$
%i.e.  $e_i %|w>=\mu_i |w>$, for some numbers $\mu_i$.
(Existence and other issues on Whittaker vectors were deeply explored
by Kostant \cite{Kost}).

The %idea to use intertwining operators to obtain
relations between Whittaker functions $W_{\lb}(\phi)$ for different $V_\lb$
will be obtained as follows: let us find the intertwinning operators $V_{\lb_1}\otimes
V_{fin} \to V_{\lb_2}$ then we will show that there exists a relations
$\cal{D}$ $(W_{\lambda_1}({\phi}_{1}, ..., \phi_{n}))= $ $ \tilde{\cal D}$
$(W_{\lambda_2} (\phi_1,..., \phi_n))$ (see formulas \ref{rais_1_sl2},\ref{rais_1_sl3}) where
${\cal D}, \tilde{\cal D}$ are some differential operators, and $V_{fin}$ is
any finite-dimensional representation of the group $SL(n+1,\RR)$.  This
relation is an example of the  so-called raising operators \cite{Rais}.  This
idea was suggested to us by S. Kharchev to whom we are deeply indebted.
Analogously we will construct bilinear relations from the intertwining
operators $V_{\lb_1}\otimes V_{\lb_2} \to V_{\lb_3} \otimes V_{\lb_4} $.  One
can proceed this way to construct other relations like above, for example
$V_{\lb_1} \otimes V_{\lb_2} = \oplus_i V_{\lb_i} $ then one can obtain
relation like $ W_{\lb_1} \otimes W_{\lb_2} = \sum_i D_i (W_{\lb_i})$
(see formula \ref{prod-sl2}).
%So one can see that the idea is quite simple and general.

Let us also note that the existence
of the bilinear identities in this situation is quite surprising,
because it means that wave-function, which are the eigen-functions
of Toda Hamiltonians, i.e. eigenfunctions of some
{\it linear} differential
operators, satisfy some system of {\it nonlinear} differential
equations. One cannot separate some finite number of wave-functions
from this system , only infinite set of them is closed nonlinear system.
Nevertheless it is some system of nonlinear equations, but we
do not know its sense, its hamiltonian description (if it is) or anything
else, it seems to be interesting to clarify it.

Let us describe what kind of technical problems arise in proceeding
to the concrete formulas.
%To write down the relation  one need the following:
Let $|w>_{\lb_1}$ be the Whittaker vector in representation $V_{\lb_1} $,
let us consider the intertwiner $\Phi : V_{\lb_1} \to V_{\lb_2} \otimes
V_{fin} $ then $ \Phi |w>_{\lb_1} = \sum_i |w>_i \otimes |i>$ where $|i>$ is
some basis in finite dimensional representation $V_{fin}$. {\it to construct
the explicit relations one needs to find such polynomials $P_i$
that  $|w>_i = P_i(f_i,h_i) |w>_{\lb_2}$ where $
P_i$ are some polynomials of generators $f_i,h_i$,
i.e.  one needs to express the components of the image of Whittaker vector in
spaces $V_{\lb_1}$ as application of some operators to Whittaker in space $
V_{\lb_2}$ }. This is the main technical problem we treat in this paper, we
solve it for the cases $SL(2,\RR)$ and $SL(3,\RR)$, but we were unable to
find the general formula for the case of $SL(n,\RR)$.

 Let us mention that the $SL(2)$ Toda equation reduces
to the well-known textbook equation, whose solutions are
Bessel and Mcdonald functions (depends on necessary asymptotics, we need
Mcdonald one), and the relation from theorem 2.11 % \ref{rais_1_sl}
one can
find, for example in \cite{Ryzik-Grads}.

 Let us argue that the scheme like above can be applied to other
integrable models. One knows that if instead of the Whittaker vector $|w>$ one
considers a spherical vector $|s>$ i.e. a vector such that it is
remains fixed under the action of the maximal compact subgroup then
one can obtain wave function for Calodgero-Sutherland integrable system -
zonal spherical function,
obviously if one considers different intertwiners one can obtain different
relations between wave functions.
Also instead of finite-dimensional Lie group $G$ one
can takes affine group $\hat{G}$, which corresponds to periodic Toda chain,
or quantum group $U_q (G)$, which leads to
difference equations \cite{Etingof99,Sevost}.
%Analogously   one can change
%group G to corresponding affine group $\tilde G$ than one obtains
%instead of open chains their periodic counterparts,  one can
%consider also quantum deformations of the group $G$ they lead
%to different difference equations \cite{Etingof99,Sevost},
Obviously,
our scheme of obtaining different relations will work in this
situations also. But of course explicit constructions of these
relations seems to be not obvious.

 We should also mention that raising operators considered
here were widely explored  recently \cite{Rais},
but all approaches are completely different from our's
and very simple idea, that the origin of this relations
is presence of the intertwining operators, was not present in the literature
to our knowledge. Let us repeat  that from our point of view
one can consider all models in the same setup.
Also we see that the  group
theory approach shows that such at first sight not related subjects such
as Baxter relation important in lattice models \cite{Feig-Anton} and
conformal field theory \cite{BLZ}, raising operators in integrable quantum
many body problems \cite{Rais} and bilinear relations in the theory
of classical nonlinear equations \cite{Kyoto,Mir} have the same explanations
from the point of view of intertwining operators.

Exposition is organized as follows: section \ref{sl(2)} is devoted to the case
of $SL(2)$ and the demonstration of the main ideas.
In the section \ref{sl(3)}
we consider the case of $SL(3)$, which is mainly analogous
to the $SL(2)$ case, but that is the case, where first
arises the problems
in writing the explicit formulas for the
action of the intertwining operators on the Whittaker vectors.
And we show how these problems can be solved for $SL(3)$.
Section \ref{sl(n)} is devoted to the $SL(N)$  case.
We recall representation theory approach to the $SL(N)$ Toda
chain and discuss the problem of obtaining explicit raising operators.
In section \ref{conclud} we make some concluding
remarks.

%%%%%%%%%%%%%%%%%%%%%%%%%%%%%%%%%%%%%%%%%%%%%%%%%%%%%%%%%%%%%%%%%%%
%%%%%%%%%%%%%%%%%%%%%%%%%%%%%%%%%%%%%%%%%%%%%%%%%%%%%%%%%%%%%%%%%%%
%%%%%%%%%%%%%%%%%%%%%%%%%%%%%%%%%%%%%%%%%%%%%%%%%%%%%%%%%%%%%%%%%%%
%%%%%%%%%%%%%%%%%%%%%%%%%%%%%%%%%%%%%%%%%%%%%%%%%%%%%%%%%%%%%%%%%%%
%%%%%%%%%%%%%%%%%%%%%%%%%%%%%%%%%%%%%%%%%%%%%%%%%%%%%%%%%%%%%%%%%%%
%%%%%%%%%%%%%%%%%%%%%%%%%%%%%%%%%%%%%%%%%%%%%%%%%%%%%%%%%%%%%%%%%%%
\section{$ SL(2)$ Toda Chain.} \label{sl(2)}
%%%%%%%%%%%%%%%%%%%%%%%%%%%%%%%%%%%%%%%%%%%%%%%%%%%%%%%%%%%%%%%%%%%
%%%%%%%%%%%%%%%%%%%%%%%%%%%%%%%%%%%%%%%%%%%%%%%%%%%%%%%%%%%%%%%%%%%
%%%%%%%%%%%%%%%%%%%%%%%%%%%%%%%%%%%%%%%%%%%%%%%%%%%%%%%%%%%%%%%%%%%
%%%%%%%%%%%%%%%%%%%%%%%%%%%%%%%%%%%%%%%%%%%%%%%%%%%%%%%%%%%%%%%%%%%
%%%%%%%%%%%%%%%%%%%%%%%%%%%%%%%%%%%%%%%%%%%%%%%%%%%%%%%%%%%%%%%%%%%
%%%%%%%%%%%%%%%%%%%%%%%%%%%%%%%%%%%%%%%%%%%%%%%%%%%%%%%%%%%%%%%%%%%

  In  this section we will show how to obtain the raising operators,
bilinear, etc. identities for the wave functions of $SL(2)$ Toda. We will
do it with the help of the intertwining operators between the tensor products
of representations of $SL(2)$. Before doing this in subsections 1 and 2
we will recall
representation theory interpretation of the Toda wave functions as matrix
element in irreducible representation, here we
closely follow \cite{Liouvil}.
Let us sketch the content: we consider Whittaker vector $|w>$ in irreducible
representation $V_{\lb}$ of $SL(2)$, i.e. such vector that it is
eigenvector for the  generator $e$, consider Whittaker function $W(\phi)=
<w| exp (\phi h) |w>$ , we show that it is eigenfunction for Toda
hamiltonian $\2 \partial_{\phi}^2 + exp(-2\phi) +\partial_{\phi}$ , the reason is
that this Hamiltonian is simply the Casimir of $SL(2)$, then we consider
intertwining operator $\Phi : V_{\lb+1} \to V_{\lb} \otimes V_{1}$,
where $V_1$ is standard two-dimensional representation of $SL(2)$,
and the main
technical point is to show that $\Phi (|w>_{\lb+1})=
\frac{1}{\lb+1}\2(h+\lb+2)|w>_{\lb} \otimes|0>+\mu_R |w>_{\lb} \otimes|1> $,
it is quite obvious for the case of $SL(2)$, because we know the representation
of the Whittaker vector in terms of the action of the generator $f$ on
the vacuum
vector in module $V_{\lb}$, the fact that such representation for the
case of $SL(N)$ is unreasonable is the main problem for writing the explicit
 formulas  for the $SL(N)$ case . Then we proceed to the bilinear relations,
to find them one should find the intertwining operator
between some representations:
$V_{\lb_1}\otimes V_{\lb_2} \to V_{\lb_3} \otimes V_{\lb_4} $.
To find such intertwiners is not an easy problem, but there is a
trick which is due to \cite{gklmm} that works as follows:
if we have intertwiners : $V_{\lb+1} \to V_{\lb} \otimes V_{1}$ and
$V_{\mu+1} \to V_{\mu} \otimes V_{1}$ we can consider
their tensor product :
$V_{\lb+1}\otimes V_{\mu+1} \to V_{\lb} \otimes V_{1} \otimes
V_{\mu} \otimes V_{1}$
then consider intertwiner: $V_{1}\otimes V_{1} \to \CC$, then
consider the composition:
$V_{\lb_1+1}\otimes V_{\mu_1+1} \to V_{\lb_1} \otimes V_{1} \otimes
V_{\mu_1} \otimes V_{1} \to V_{\lb+1} \otimes V_{\mu+1}$, hence we get
desirable intertwiner (consideration instead of $V_{1} $ other
finite dimensional representations leads to various kinds of necessary
intertwiners).

\subsection {Notations:}

Let us denote by $ (\lb,n)$ the product $(\lb)(\lb-1)(\lb-2)...(\lb-n+1)$,
we will often use this abbreviation.

The Lie algebra $SL(2)$ is defined by the generators $e,f,h$
satisfying the following relations:

$$
[e,f]=h,\ \ [h,e]=2e, \ \ [h,f]=-2f.\\
$$

The fundamental representation:
%$$ %

\bea
h=\pmatrix 1&0\\0&-1 \endpmatrix,\
e=\pmatrix 0&1\\0& 0 \endpmatrix,\
f=\pmatrix 0&0\\1& 0 \endpmatrix.
\eea
%\bea
%h=\left(\begin{array}{cc}1&0\\0&-1\end{array}\right),\
%e=\left(\begin{array}{cc}0&1\\0&0\end{array}\right),\
%f=\left(\begin{array}{cc}0&0\\1&0\end{array}\right).
%$$ %
%\eea

The quadratic Casimir operator:
\bea
C=(ef+fe)+\2 h^2=2fe+h+\2 h^2.
\eea

{Highest weight representation:}

We consider the principal (spherical) series of
representations, induced by the
one-dimensional representations of the Borel subgroup.
The space of representation $V_\lb$ is functions of one real variable $x$
and matrix elements are defined by integrals with the flat measure.
The action of the algebra is given by differential operators:
\bea\label{bwb-sl2}
e = \frac{\partial}{\partial x},\ h=-2x{\d x}+\lb,\ f=-x^2{\d x}+\lb x.
\eea
The action of the group is given by the formula:
\bea
\pmatrix a & b \\ c & d \endpmatrix : x^k \to  (a+cx)^{\lb-k}(b+dx)^k
\eea

Obviously that $Constants$ are vacuum (i.e. the highest weight)
vectors for such representation i.e.
$e (Constant) = 0,~~~~h (Constant)=\lb (Constant)$

\subsection{Whittaker function.}
%%%%%%%%%%%%%%%%%%%%%%%%%%%%%%%%%%%%%%%%%%%%%%%%%%%%%%%%%%%%%%%%%%%
%%%%%%%%%%%%%%%%%%%%%%%%%%%%%%%%%%%%%%%%%%%%%%%%%%%%%%%%%%%%%%%%%%%

{\lD 2.1:}
The vector $|w>_{\lb}^{\mu}$ is called {\it Whittaker vector} in representation $V_\lb$,
if it is eigenvector for generator $e$ i.e.
$e |w>_{\lb}^{\mu}=\mu |w>_{\lb}^{\mu}$ for some
constant $\mu$.

{\lL 2.2:}
For given $\mu$ and irreducible representation $V_\lb$~ ($\lb \ne 1,2,3,...$)
there exists a unique (up to a scalar)
Whittaker vector $|w>_{\lb}^{\mu} \in V_\lb$, which can be expressed as follows:
\bea
|w>_{\lb}^{\mu} = exp (\mu x) = \sum_{i=0} \frac{\mu^n f^n}{n! (\lb,n)} |vac>
\label{wh_sl2}
\eea

where we denote $ (\lb,n)=\lb (\lb-1) ... (\lb-n+1)$ and $|vac>$ - vacuum
( the highest weight) vector in $V_\lb$.

The proof is obvious.

The dual Whittaker vector: $~~{}_{\lb}^{\mu}<w|$ is given
$\sum_{n=0}^{\infty} {}_{\lb}<0|\frac{\mu^n e^n}{n! (\lb,n)} = x^{-\lb-2}
exp(-\mu x)$
 (see \cite{Liouvil}).

{\lD 2.3:} The Whittaker function $W_{\lb}^{\mu_L,\mu_R}(\phi)$ is the function given by:
\bea
W_\lb^{\mu_L,\mu_R} (\phi)= ~~{}_{\lb}^{\mu_L}<w|exp( \phi h)|w>_{\lb}^{\mu_R}
\label{wf_sl2}
\eea

{\Rem} we will sometimes omit indexes $\lb, \mu$ of Whittaker vectors and
Whittaker functions, if they are not important.

{\lPn 2.4:}
$
W_{\lb}(\phi)=
exp(\lb \phi) \sum_{i=0}^{\infty}
\frac{(\mu_R)^i (\mu_L)^n  exp(-2i\phi)}{i!(\lb,i)}
$

The proof is obvious recalling (\ref{wf_sl2}) and $<vac|e^n|f^n|vac>=n!(\lb,n)$

\Rem  for the case $SL(2,\RR)$ Whittaker function
coincides with the textbook Macdonald function (up to some exponents and scalars),
this can be easy seen from the integral representation below, see also \cite{Liouvil}.

%{\it Remark:} this notation is not accordered with Lang's "$\bf{SL(2,\RR)}$",
%but follows the modern notation (for example \cite{Liouvil},\cite{Etingof99}).

{\lPn 2.5:}
The Whittaker function satisfies the equation:
\bea\label{HI2}
\left[\2 {\partial^2\over\partial\phi^2}+{\d\phi}-2\mu_R\mu_L e^{-2\phi}\right]
W_{\lb}(\phi)=  (\2\lb^2 +\lb) W_{\lb}(\phi)
\eea

Hence, the function $\Psi_{\lb} (\phi)=exp(\phi)W_{\lb}(\phi)$ satisfies
two  particle open Toda equation:
\bea\label{SchG}
\left[\2
{\partial^2\over\partial\phi^2}-2\mu_R\mu_Le^{-2\phi}\right]
\Psi_{\lb}(\phi)=\left(\2\lb^2+\lb\right)\Psi_{\lb} (\phi),
\eea

The proof is very simple, but we will reproduce it, since it is the same
for all kinds of integrable many body problems.
The main idea is that the Hamiltonian is a Casimir operators, hence
it acts by scalar on every matrix element, that's how RHS of \ref{HI2}
appears, on the other hand it is some concrete operator - that's
how LHS of \ref{HI2} appears.

{\lP}
The Casimir operator $C= 2fe+h+\2 h^2$ acts as scalar $\2\lb^2+\lb$
on the representation $\vL$. Hence:
$
<w|exp (\phi h) C \w= (\2\lb^2+\lb) W_{\lb}(\phi)
$
On the other hand   $ <w|exp (\phi h) C \w=
<w|exp (\phi h) (\2 h^2+h+2fe) \w=
(\2 {\partial^2\over\partial\phi^2}+{\partial\over\partial\phi})
W_{\lb}(\phi)+ 2<w|f~ exp (\phi(h-2)) e\w=
(\2 {\partial^2\over\partial\phi^2}+{\partial\over\partial\phi}
+ 2\mu_R \mu_l exp (-2\phi) )W_{\lb}(\phi)
$ hence we obtain RHS of \ref{HI2}.
{\QD}

\Rem any matrix element of the irreducible representation
is eigenfunction if all the Casimir operators of any group.
It may seem that hence we obtain integrable system. Of course, it's not
true,  because usually the dimension of a group is something
like $n^2$ where the number of independent Casimirs is something
like $n$. So one does not have enough Hamiltonians in involution.
So one should reduce the phase space from the hole
group to some submanifold of it, but  the problem is that
Casimir operator acting on the function on this submanifold
can give a function beyond that class.
In our case we saw that $<w|exp (\phi h) C \w=
(\2 {\partial^2\over\partial\phi^2}+{\partial\over\partial\phi}
+ 2\mu_R \mu_l exp (-2\phi) )W_{\lb}(\phi)
$. This is due to the clever choice: we considered
matrix element of the Cartan group element
between Whittaker vectors, the idea to do so is due
to   Kazhdan and Kostant.
For example if one takes some arbitrary vectors $<u|, |v>$
then
 $<u|exp (\phi h) C |v>=
(\2 {\partial^2\over\partial\phi^2}+{\partial\over\partial\phi})W_{\lb}(\phi)
+ 2<u|exp (\phi h) fe |v>
$, and one can do nothing with the  $<u|exp (\phi h) fe |v>$
one cannot in general express it via $<u|exp (\phi h) |v>$.
So one cannot obtain the LHS of (\ref{HI2}), despite one has RHS.
The general conception
of the consistent hamiltonian reductions is due to Drinfeld and Sokolov \cite{Dr-Sok}.

{\lPn 2.6 (Integral representation of the Whittaker function) :}
\bea \label{LWFG}
W_\lb(\phi)=
\left({1\over \mu_R}\right)^{-(\lb+1)}
\int_0^{\infty} x^{-(\lb+2)}e^{-{\mu_L\mu_Re^{-2\phi}
\over x}- x}dx=%\\
2 e^{(\lb+1)\phi}\left(\sqrt{{\mu_L\over\mu_R}}\right)^{-(\lb+1)}
K_{\lb+1}(2\sqrt{\mu_L\mu_R}e^{-\phi})\ dx.
\eea

- where $K_{\lb} (z)$ is the Macdonald function.

\Rem for the convergence we should require $Re \mu_L\mu_R>0$.

The theorem is tautology after recalling that the invariant pairing
$<v|u>$ is realized as integral $\int_0^{\infty}$, and
dual Whittaker vector $<w|$ is realized as $x^{-2(\lb+1)}e^{-{\mu_L\over
x}}$.  The second equality is the textbook integral formula for the Macdonald
function.

\subsection{Intertwining operators}
%%%%%%%%%%%%%%%%%%%%%%%%%%%%%%%%%%%%%%%%%%%%%%%%%%%%%%%%%%%%%%%%%%%
%%%%%%%%%%%%%%%%%%%%%%%%%%%%%%%%%%%%%%%%%%%%%%%%%%%%%%%%%%%%%%%%%%%

Denote by $V_1$ the fundamental representation of the $SL(2,\RR)$,
denote  $|0>$ the highest weight vector, $|1>$ the lowest weight vector, i.e.
\bea
|0>=\left(\begin{array}{cc}1\\0\end{array}\right),\
|1>=\left(\begin{array}{cc}0\\1\end{array}\right),\
%$$ %
\eea
hence $e |0>=0, f|0>=|1>, h |0>=|0>; e |1>=|0>, f|1>=0, h |1>=-|0>$.

%Let us construct the  intertwining operators $V_{\lb+1} \to \vL\otimes V_1$

{\lPn 2.7:}
The isomorphism $\Phi_{\lb} = \Phi_{\lb,+} \oplus \Phi_{\lb,-} :
V_{\lb+1} \oplus V_{\lb-1} \to \vL\otimes V_1   $ is given
 by the formulas:
\bea
\label{int_sl2}
&&\Phi_{\lb,+} : f^n|vac>_{\lb+1} \to f^n|vac>_{\lb}\otimes |0> +
n f^{n-1}|vac>_{\lb}\otimes |1>,\\
&&
\Phi_{\lb,+}^{-1} : f^n|vac>_{\lb}\otimes |0> \to \frac{(\lb+1-n)}{\lb+1} f^n |vac>_{\lb+1},
~~~~~~~~~\Phi_{\lb,+}^{-1} : f^n|vac>_{\lb}\otimes |1> \to \frac{1}{\lb+1}f^{n+1} |vac>_{\lb+1},\\
&&
\Phi_{\lb,-} : f^n|vac>_{\lb-1} \to f^{n+1}|vac>_{\lb}\otimes |0> +
(n-\lb) f^{n}|vac>_{\lb}\otimes |1>,\\
&&
\Phi_{\lb,-}^{-1} : f^n|vac>_{\lb}\otimes |0> \to \frac{n}{\lb+1} f^{n-1} |vac>_{\lb-1},
~~~~~~~~~\Phi_{\lb,-}^{-1} : f^n|vac>_{\lb}\otimes |1> \to \frac{-1}{\lb+1} f^{n} |vac>_{\lb-1}.
\eea
The action of dual operators $\Phi^{*},\Phi^{*}_{1},\Phi^{*}_{-1}$ on the
dual anti-representations are given by the same formulas (up to scalars) with
the change of $f$ to $ e$ , which follows from the Chevalley antiinvolution:
\bea
\label{int*_sl2}
&&\Phi_1^* : {}_{\lb+1}<vac|e^n \to
{}_{\lb}<vac|e^n \otimes <0| +
n{}_{\lb}<vac|e^{n-1} \otimes <1|\\
&&
{\Phi_1^{-1}}^* : {}_{\lb}<vac|e^n\otimes <0| \to
\frac{(\lb+1-n)}{\lb+1}~~ {}_{\lb+1}<vac|e^n ,
~~~~~~~~~{\Phi_1^{-1}}^* : {}_{\lb}<vac|e^n\otimes <1| \to
\frac{1}{\lb+1}~~ {}_{\lb+1}<vac|e^{n+1} ,\\
&&
\Phi_{-1}^* : {}_{\lb-1}<vac|e^n \to
\frac{1}{\lb+1} ({}_{\lb}<vac|e^{n+1} \otimes <0| +
(n-\lb){}_{\lb}<vac|e^{n} \otimes <1|)\\
&&
{\Phi_{-1}^{-1}}^* : {}_{\lb}<vac|e^n\otimes <0| \to
(n \lb )~~{}_{\lb-1}<vac|e^{n-1} ,
~~~~~~~~~{\Phi_{-1}^{-1}}^* : {}_{\lb}<vac|e^n\otimes <1| \to
(-\lb) ~{}_{\lb-1}<vac|e^{n} .
\eea
The invariant pairing is given by: $<0|0>=1, <0|1>=0, <1|0>=0,
<1|1>=1.$

The proof is obvious, let us only note that in order
to find the intertwiners $\Phi_{\pm1}^{-1}$ from $V_{\lb \pm 1}\otimes V_1
\to V_{\lb}$ it is more convinient to express them via the dual intertwiners:
$V_{\lb \pm 1} \to  V_1^{*} \otimes  V_{\lb}$, which can be easily found.

{\lPn 2.8:}
Intertwiners act on Whittaker vectors as follows:
\bea
%&&
%\Phi^{-1}_{-1} : |w>_{\lb}\otimes|0>\to \frac{\mu}{\lb} |w>_{\lb-1},
%~~~~~~
&&
%\Phi_1 :
|w>_{\lb+1}\to \frac{(h+\lb+2)}{2(\lb+1)}|w>\otimes|0>+
\frac{\mu_R}{\lb+1} |w>\otimes|1>,\\
&&
%\Phi_{-1} :
|w>_{\lb-1}\to \frac{\lb}{\mu_R}\2(-h+\lb)|w>\otimes|0>
-\lb |w>\otimes|1>,\\
&&
|w>_{\lb}\otimes |0> \to |w>_{\lb+1}+ \frac{\mu_R}{\lb(\lb+1)}|w>_{\lb-1},\\
&&
|w>_{\lb}\otimes |1> \to \frac{(\lb+1-h)}{2\mu_R} |w>_{\lb+1}-
\frac{\lb+1+h}{2\lb}|w>_{\lb-1}, \\
&&
%{\Phi^{-1}}^{*} :
~ {}_{\lb}<w|\otimes <0|\to
\mu_L~~ {}_{\lb-1}<w|+
~~{}_{\lb+1}<w|,\\
&&%{\Phi^{-1}}^{*} :
~~{}_{\lb}<w|\otimes <1|\to
-\frac{h+\lb+1}{2\lb}~~{}_{\lb-1}<w|\otimes <1| +
\frac{\lb+1-h}{2\mu_L}~~ {}_{\lb+1}<w|.
\eea
\lP
The proof is simple, let us consider  only the first
equality.
We will use the following  trivial fact:

$(\lb-n+1) f^n |vac>_{\lb}=
\2(\lb+2+h) f^n |vac>_{\lb}$.

\noindent
$\Phi_{1}(|w>_{\lb+1})=\Phi_{1}(\sum_{n=0}^{\infty}\frac{\mu_R^n}{n!(\lb+1,n)}
f^n|vac>_{\lb+1})=
\sum_{n=0}^{\infty}\frac{\mu_R^n}{n!(\lb+1,n)}(
f^{n}|vac>_{\lb}\otimes |0>+nf^{n-1}|vac>_{\lb}\otimes |1>)=$

$=
\frac{1}{\lb+1}
\sum_{n=0}^{\infty}\frac{\mu_R^{n}}{(n)!(\lb,n)}
(\lb-n+1) f^{n}|vac>_{\lb}\otimes |0>
+\frac{\mu_R}{\lb+1} |w>_{\lb}\otimes|1>=
\frac{(\lb+2+h)}{2(\lb+1)} |w>_{\lb}\otimes|0>+
\frac{\mu_R}{\lb+1} |w>_{\lb}\otimes|1>.
$

~~\QD

\Rem let us note that such simple formulas like above cannot be true for an arbitrary
vector $|v>_{\lb}\in \vL$ i.e. of course  from the irredicibility of repsentations follows
that there exist series $O(f,h,e)$ such that $\Phi_1(|v>_{\lb+1})
=O(f,h,e) |v>_{\lb}$,
but in general such expression will be complicated.

Let us find the intertwiners between tensor product of
infinite-dimensional representations.

\lPn 2.9: the isomorphism $\Phi=\oplus_{k=0}^{\infty} \Phi_k: \vL\otimes
V_{\nu} \to \oplus_{k=0}^{\infty} V_{\lb+\nu-k}$ in realization (\ref{bwb-sl2})
is given by the formulas:
\bea
\Phi_k : x^n \otimes x^m \to
C_k \sum_{i=0}^{k} (-1)^i
\left(\begin{array}{cc}n\\k-i\end{array}\right)
\left(\begin{array}{cc}m\\i\end{array}\right)
\frac{(\nu)_k}{(\nu)_i}\frac{(\lb)_k}{(\lb)_{k-i}} x^{n+m-k}, \\
where~~~~~ C_k=\frac{ \left(
\prod_{i=1}^{k} (i (\nu^2(\lb-i+1)+\lb^2(\nu-i+1)) \right)^{\2} }
{ (\nu \lb)^k (\sum_{i=0}^{k}
\frac{(\nu)_k}{(\nu)_i}\frac{(\lb)_k}{(\lb)_{k-i}} ) }.
\eea
Or one can rewrite it as follows. Denote by $mul_k$ operator
$\vL\otimes V_{\nu} \to  V_{\lb+\nu-k}$  which acts:
$f(x)\otimes g(x)\to f(x)g(x)$. Then:
\bea
\label{int-sl2-2-2}
\Phi_k: f(x)\otimes g(x) \to  C_k(mul_k)*(
\sum_{i=0}^k \frac{ (-1)^i}{(k-i)! i!}
\frac{(\nu)_k}{(\nu)_i}\frac{(\lb)_k}{(\lb)_{k-i}}
\partial_{x}^{k-i} f(x) \otimes \partial_{x}^{i} g(x) ).
\eea

The proof can be found in the literature,  or proceeded as follows:
from the equality: $\Phi_k* h=h*\Phi_k$ follows
$\Phi_k(x^n\otimes x^m)=c_{n,m}x^{n+m-k}$, using $\Phi_k* e=e*\Phi_k$
 one obtains the recurrence relation for the coefficients $c_{n,m}$.
 Solution of this relation can be guessed after considering several examples
of small $k=0,1,2$.

\lPn 2.10: The action of $\Phi_k$ on Whittaker vectors
is given by:
\bea
\Phi_k: |w>_{\lb}^{\mu_1}\otimes |w>_{\nu}^{\mu_2} \to  C_k
\sum_{i=0}^k \frac{ (-1)^i}{(k-i)! i!}
\frac{(\nu)_k}{(\nu)_i}\frac{(\lb)_k}{(\lb)_{k-i}}
{(\mu_1)}^{k-i} {(\mu_2)}^{i})
|w>_{\lb+\nu-k}^{\mu_1+\mu_2} .
\eea
This easily follows from \ref{int-sl2-2-2}.

\subsection{Relations for the Whittaker wave functions.}
%%%%%%%%%%%%%%%%%%%%%%%%%%%%%%%%%%%%%%%%%%%%%%%%%%%%%%%%%%%%%%%%%%%
%%%%%%%%%%%%%%%%%%%%%%%%%%%%%%%%%%%%%%%%%%%%%%%%%%%%%%%%%%%%%%%%%%%

{$~ $ }

{\lT (Raising operators) 2.11} The following relations holds:
\bea
\label{rais_1_sl2}
&&
W_{\lb+1}(\phi)=
exp(\phi)\frac{\partial_{\phi}+\lb+2}{2(\lb+1)} W_{\lb}(\phi),\
\mu_L W_{\lb-1}(\phi)= \frac{\lb}{\mu_R}
exp(\phi)\frac{-\partial_{\phi}+\lb}{2} W_{\lb}(\phi).
\eea
\Rem easy to see that these relations are consistent: application of one
then another gives that $W_{\lb}$ satisfy Toda equation.

\lP
We will calculate $ {}_{\lb}<w|\otimes <0| ~~~\Phi_{1} exp(\phi h) |w>_{\lb+1}$
in  two ways, first we will apply $\Phi_{1}$ to $exp(\phi h) |w>_{\lb+1}$
and then take pairing, hence we will obtain the RHS of (\ref{rais_1_sl2});
to obtain the LHS of (\ref{rais_1_sl2}) we will apply $\Phi_{1}^{*}$ to
$ {}_{\lb}<w| \otimes <0| $ and consider the pairing after.

$ {}_{\lb}<w|\otimes <0| ~~~\Phi_{1} exp(\phi h) |w>_{\lb+1}=
{}_{\lb}<w|\otimes <0| ~~~exp(\phi h) \Phi_{1}  |w>_{\lb+1}=
{}_{\lb}<w|\otimes <0| ~~~exp(\phi h)
(\frac{(h+\lb+2)}{2(\lb+1)}|w>_{\lb}\otimes|0>+\frac{\mu_R}{\lb+1} |w>\otimes|1>)=
{}_{\lb}<w|\otimes <0| ~~~
\frac{(h+\lb+2)}{2(\lb+1)} exp(\phi h)|w>_{\lb} \otimes exp(\phi h)|0>=
exp(\phi) \frac{(\partial_{\phi}+\lb+2)}{2(\lb+1)}
{}_{\lb}<w|\otimes <0| ~~~exp(\phi h)|w>_{\lb} .$
Hence we get RHS of \ref{rais_1_sl2}.

$ {}_{\lb}<w|\otimes<0| ~~~\Phi_{1} exp(\phi h) |w>_{\lb+1}=
 \left( {}_{\lb}<w|\otimes<0| \Phi_{1}^{*}\right) | exp(\phi h) |w>_{\lb+1}=
 {}_{\lb+1}<w| exp(\phi h) |w>_{\lb+1}= W_{\lb+1}(\phi).$
Hence we get the LHS of  (\ref{rais_1_sl2}). \QD The proof of the second equality in
(\ref{rais_1_sl2}) is the same.

One can deduce directly
from the (\ref{rais_1_sl2}) or deduce
representation theoretically (as we will do for the demonstration
of the idea) the following  relations:

\lT (Baxter like relations) 2.12:
\bea
&& exp(\phi) W_{\lb}(\phi)= \frac{\mu_L \mu_R }{\lb(\lb+1)}
W_{\lb-1}(\phi)  +  W_{\lb+1}(\phi) ,\nn\\
&&exp(-\phi) W_{\lb}(\phi)=
\left( \frac{\partial_{\phi}+\lb+1}{2} \right)^2 W_{\lb-1}(\phi)+
 \frac{1}{\mu_L \mu_R} \left(\frac{\lb+1-\partial_{\phi}}{2} \right)^2
W_{\lb+1}(\phi).
\eea
\Rem analogous considerations for the infinite-dimensional algebras
leads to Baxter relations for the quantum transfer matrices in lattice
models, see \cite{Feig-Anton}.

\lP
We will use isomorphism $\Phi: V_{\lb+1}\oplus
V_{\lb-1} \to V_{\lb}\otimes
V_1$ and calculate the pairing
$<0|<w| exp (\phi h) |w>|0>$ first directly, obtaining LHS, then applying
isomorphism $\Phi^{-1}$ to $ |w>|0>$ and ${\Phi^{-1}}^{*}$ to $ <w|<0|$, and
obtain RHS.
\bea
&&exp(\phi) W_{\lb}(\phi)=
 {}_{\lb}<w|\otimes<0| exp(\phi h)|0>\otimes|w>_{\lb}=\nn\\
&&
(\mu_L~~ {}_{\lb-1}<w|+ ~{}_{\lb+1}<w|)
exp(\phi h) (\frac{\mu_R}{\lb(\lb+1)}|w>_{\lb-1}+|w>_{\lb+1})=\nn \\
&&
=\frac{\mu_L \mu_R}{\lb(\lb+1)}~ {}_{\lb-1}<w| exp(\phi h)|w>_{\lb-1}+
~ {}_{\lb+1}<w| exp(\phi h)|w>_{\lb+1}=
\frac{\mu_L \mu_R}{\lb(\lb+1)} W_{\lb-1}(\phi)+W_{\lb+1}(\phi).
\eea
\QD
The proof of the second relation is the same, but one should consider
${}_{\lb}<w|\otimes<1| exp(\phi h)|1>\otimes|w>_{\lb}$ instead of
${}_{\lb}<w|\otimes<0| exp(\phi h)|0>\otimes|w>_{\lb}$. Obviously
$0={}_{\lb}<w|\otimes<0| exp(\phi h)|1>\otimes|w>_{\lb}$, so one cannot
obtain any more relations like above.

\lT (Bilinear relations) 2.13:
\bea
\label{bil_sl2}
&&\left(
-\frac{\mu^R_2 exp(\phi_1-\phi_2)}{2(\lb+1)(\nu+1)}(\partial_{\phi_1}+\lb+2)+
\frac{exp(\phi_2-\phi_1)\mu^R_{1}}{2(\lb+1)(\nu+1)}
(\partial_{\phi_2}+\nu+2)
\right) W_{\lb}^{\mu^L_1,\mu^R_1}(\phi_1) W_{\nu}^{\mu^L_2,\mu^R_2}(\phi_2)
= \\
&&=
\left(- \frac{\nu+1-\partial_{\phi_2}}{2\mu_2^L}+
\frac{\lb+1-\partial_{\phi_1}}{2\mu^L_1}
\right)  W_{\lb+1}^{\mu^L_1,\mu^R_1}(\phi_1)
W_{\nu+1}^{\mu^L_2,\mu^R_2}(\phi_2)
\eea

\lP
let us denote by $\Phi_{\lb,\nu}$ the intertwiner:
 $ V_{\lb+1}\otimes V_{\nu+1} \to V_{\lb}\otimes V_{\nu} $.
we will construct that intertwiner using a trick due to \cite{gklmm}:
let us denote by $S$ the intertwiner $V_1\otimes V_1 \to \CC$,
than $\Phi_{\lb,\nu}$ is given by composition:
$ (id\otimes S \otimes id) (\Phi_{\lb,+}\otimes \Phi_{\nu,+}) $.

\lL 2.14: The action of the $\Phi_{\lb,\nu}$ on Whittaker vectors
is given by the formulas:
\bea
&&
\Phi_{\lb,\nu} : exp (\phi_1 h) \w_{\lb+1}^{\mu_1^R} \otimes
exp (\phi_2 h) \w_{\nu+1}^{\mu_2^R} \to \nn\\
&&
-exp(\phi_1 - \phi_2) \frac{ h+\lb+2}{2(\lb+1)}exp(\phi_1 h)\wl^{\mu_1^R}\otimes
\frac{ \mu^R_2 exp(\phi_2 h)}{\nu+1}exp(\phi_2 h)\w_{\nu}^{\mu_2^R}+\\
&&
+exp(\phi_2 - \phi_1)  \frac{\mu^R_1 exp(\phi_1 h)}{(\lb+1)}
\wl^{\mu_1^R}\otimes \frac{(h+\nu+2)}{2(\nu+1)} exp(\phi_2 h)\w_{\nu}^{\mu_2^R},\\
&&
\Phi_{\lb,\nu}^{*}:~~
{}^{\mu_1^L}\wwl\otimes~~ {}_{\nu}^{\mu_2^L} \ww \to
- ~~{}_{\lb+1}^{\mu_1^L}\ww \otimes  ~~ {}_{\nu+1}^{\mu_2^L}\ww \frac{\nu+1-h}{2\mu^L_2} +
~~ {}_{\lb+1}^{\mu_1^L}\ww \frac{\lb+1-h}{2\mu^L_1}\otimes~~ {}_{\nu+1}^{\mu_2^L}\ww.
\eea
The proof of the lemma follows from  (\ref{int_sl2}) and recalling the
fact that action of $S$  is given by:
$ S(|0>|1>)=0,~~~S(|1>|0>)=0,~~~S(|1>|0>)=1,~~~S(|0>|1>)=-1$ and respectively:
$ S^*: \CC\to V\otimes V$ acts as follows: $ 1 S^*=-<0|<1|+<1|<0|$.

The proof of the theorem easily follows from the lemma:
$
\wwl \otimes {}_{\nu}\ww \left(\Phi_{\lb,\nu} exp (\phi_1 h) \w_{\lb+1} \otimes
exp (\phi_2 h) \w_{\nu+1} \right)$ equals to LHS of \ref{bil_sl2} and
$\left( \wwl \otimes {}_{\nu}\ww  \Phi_{\lb, \nu}^{*} \right)
exp (\phi_1 h) \w_{\lb+1} \otimes
exp (\phi_2 h) \w_{\nu+1} $ equals to RHS of \ref{bil_sl2}.

And obviously
$
\wwl \otimes {}_{\nu}\ww \left(\Phi_{\lb,\nu} exp (\phi_1 h) \w_{\lb+1} \otimes
exp (\phi_2 h) \w_{\nu+1} \right)=
\left( \wwl \otimes {}_{\nu}\ww  \Phi_{\lb, \nu}^{*} \right)
exp (\phi_1 h) \w_{\lb+1} \otimes
exp (\phi_2 h) \w_{\nu+1} $.
\QD

\lT 2.15 (Nonlinear equations):
\bea
\label{nonl_sl2}
&&
\frac{(\partial_{\phi}+\lb+2)}{2(\lb+1)} W_{\lb}(\phi)
\frac{\mu_2^R(1- \partial_{\phi})}{\nu+1} W_{\nu}(\phi)
+
\frac{\mu^R_1}{\lb+1} W_{\lb}(\phi)
\frac{(\partial_{\phi}+\lb+2)}{\nu+1}(1+\partial_{\phi})
 W_{\nu}(\phi)=\\
&&
=
-W_{\lb+1}(\phi)
\frac{\nu+1-\partial_{\phi}}{2\mu^L_1}\partial_\phi W_{\nu+1}(\phi)
+
\frac{\lb+1-\partial_{\phi}}{2\mu^L_1}W_{\lb+1}(\phi)
\partial_\phi W_{\nu+1}(\phi).
\eea

The proof is the standard procedure of deducing nonlinear equations
from the bilinear  ones: one should substitute $\phi_2=\phi_1+\delta$
and consider the Taylor expansion of the $W_(\phi_2)$ in $\delta$ and
necessary relations are obtained as equalities in different powers
$\delta^k$, the above one is for $k=1$.
In this way one obtains Kadomtsev-Petviashvili hierarchy
from the Hirota bilinear relations.

\lT 2.16 (Product formula):
\bea
\label{prod-sl2}
&&
W_{\lb}^{\mu_L,\mu_R} (\phi) W_{\nu}^{\mu_L^{\prime},\mu_R^{\prime}} (\phi)=
\\
&&=
\sum_{k=0}^{\infty}
C_k W_{\lb+\nu-k}^{\mu_L+\mu_L^{\prime}, \mu_R+\mu_R^{\prime}}(\phi)
\sum_{i=0}^{k} \frac{ (-1)^i\mu_R^i{(\mu_R^{\prime})}^{(k-i)}}{(k-i)!i!}
\left(\begin{array}{cc}n\\k-i\end{array}\right)
\left(\begin{array}{cc}m\\i\end{array}\right)
\sum_{i=0}^{k} \frac{ (-1)^i\mu_L^i{(\mu_L^{\prime})}^{(k-i)}}{(k-i)!i!}
\left(\begin{array}{cc}n\\k-i\end{array}\right)
\left(\begin{array}{cc}m\\i\end{array}\right) \nn
\eea
To prove the theorem above one should note that LHS in \ref{prod-sl2}
equals to
${}_{\lb}<w|\otimes {}_{\nu}<w| exp(\phi h) \wl \otimes \w_{\nu}$, on the
other hand one can compute this expression by use of isomorphism
from proposition 2.9 and obtain RHS of  \ref{prod-sl2}.
The calculation is very similar to
the one in proof of the theorem  2.12 (Baxter relation).

So in this section we demonstrated how to deduce different
relations for the Whittaker wave function of Toda by use of
intertwiners.

%%%%%%%%%%%%%%%%%%%%%%%%%%%%%%%%%%%%%%%%%%%%%%%%%%%%%%%%%%%%%%%%%%%
%%%%%%%%%%%%%%%%%%%%%%%%%%%%%%%%%%%%%%%%%%%%%%%%%%%%%%%%%%%%%%%%%%%
%%%%%%%%%%%%%%%%%%%%%%%%%%%%%%%%%%%%%%%%%%%%%%%%%%%%%%%%%%%%%%%%%%%
%%%%%%%%%%%%%%%%%%%%%%%%%%%%%%%%%%%%%%%%%%%%%%%%%%%%%%%%%%%%%%%%%%%
%%%%%%%%%%%%%%%%%%%%%%%%%%%%%%%%%%%%%%%%%%%%%%%%%%%%%%%%%%%%%%%%%%%
%%%%%%%%%%%%%%%%%%%%%%%%%%%%%%%%%%%%%%%%%%%%%%%%%%%%%%%%%%%%%%%%%%%
\section{$ SL(3)$ Toda Chain.} \label{sl(3)}
%%%%%%%%%%%%%%%%%%%%%%%%%%%%%%%%%%%%%%%%%%%%%%%%%%%%%%%%%%%%%%%%%%%
%%%%%%%%%%%%%%%%%%%%%%%%%%%%%%%%%%%%%%%%%%%%%%%%%%%%%%%%%%%%%%%%%%%
%%%%%%%%%%%%%%%%%%%%%%%%%%%%%%%%%%%%%%%%%%%%%%%%%%%%%%%%%%%%%%%%%%%
%%%%%%%%%%%%%%%%%%%%%%%%%%%%%%%%%%%%%%%%%%%%%%%%%%%%%%%%%%%%%%%%%%%
%%%%%%%%%%%%%%%%%%%%%%%%%%%%%%%%%%%%%%%%%%%%%%%%%%%%%%%%%%%%%%%%%%%
%%%%%%%%%%%%%%%%%%%%%%%%%%%%%%%%%%%%%%%%%%%%%%%%%%%%%%%%%%%%%%%%%%%

  In  this section we will show how to obtain raising operators and
the bilinear identities for the wave functions of $SL(3)$ Toda.
The idea is completely analogous to the case of $SL(2)$:
the wave functions are matrix elements of $SL(3)$ in irreducible
representations,  so considerations of the different intertwiners:
$\vL\otimes V_{finite-dim}\to V_{\lb\pm 1}$,
$V_{\lb_1} \otimes V_{\lb_2} \to V_{\nu_1} \otimes V_{\nu_2} $, etc.
gives different relations for the wave functions.
The main technical problem here
is to find the explicit
%realization of such relations
%one should
expression for the image of the Whittaker vector $\wl$
through the Whittaker vector $\w_{\lb\pm1}$: find $P(f,h)$,
such that $\Phi (\wl)=P(f,h)(\w_{\lb+1}$, pay attention that we need
expression of the operator  $P$ in terms of generators of the
$SL(3)$, (expression in operators $x_{i,j}, \partial_{x_{i,j} }$ in
Borel-Weil realization one can easily find).
In the case of
$SL(2)$ such problem was easily solved by the fact that
there was an explicit  realization of the
vector $\wl$ :
 $\wl=\sum_{i=0} \frac{\mu^n f^n}{n! (\lb,n)} |vac>_\lb$,
but  in case of $SL(3)$ such formula is not reasonable or
useable. So one need to find the main formula:
$\Phi (\wl)=P(f,h)(\w_{\lb+1}$ using more or less the only information about
vector $\wl$ that it is eigenvector of the generators $e_i$.
So it is the main difference with the case of $Sl(2)$ and we will
mostly pay attention to this question and we will be brief in questions
analogous to the case of $SL(2)$.

\subsection {Notations:}
%%%%%%%%%%%%%%%%%%%%%%%%%%%%%%%%%%%%%%%%%%%%%%%%%%%%%%%%%%%%%%%%%%%
%%%%%%%%%%%%%%%%%%%%%%%%%%%%%%%%%%%%%%%%%%%%%%%%%%%%%%%%%%%%%%%%%%%

%Let us denote by $ (\lb,n)$ the product $(\lb)(\lb-1)(\lb-2)...(\lb-n+1)$.

The Lie algebra $SL(3)$ is defined by the generators
$e_1,e_2,e_{12},h_1,h_2,f_1,f_2,f_{12}$. Sometimes we will refer
$e_{12}$ as $e_3$, and $f_{12}$ as $f_3$. Generators
satisfy the following relations:

\bea
[e_1,f_1]=h_1,\ \ [h_1,e_1]=2e_1, \ \ [h_1,f_1]=-2f_1,\\ ~
[e_2,f_2] = h_2,\ \ [h_2,e_2]=2e_2, \ \ [h_2,f_2]=-2f_2,\\ ~
[e_{12},f_{12}] = h_1+h_2,\ \ [h_1+h_2,e_{12}]=2e_{12}, \ \ [h_1+h_2,f_{12}]=-2f_{12},\\ ~
[h_1,e_2] = -e_2,\ \ [h_1,f_2]=f_2,\ \  [h_2,f_1]=f_1,\ \   [h_2,e_1]=e_1,\\ ~
[e_1,e_2] = e_{12},\ \ [e_1,e_{12}]=0,\ \   [e_2,e_{12}]=0,\\ ~
-[f_1,f_2] = f_{12},\ \ [f_1,f_{12}]=0,\ \ [f_2,f_{12}]=0,\\ ~
[e_1,f_2]=0,\ \ [e_2,f_1]=0.\\
\eea

The first fundamental representation:
\bea
h_{1}=\left(\begin{array}{ccc}1&0&0\\0&-1&0\\0&0&0\end{array}\right),\
h_2=\left(\begin{array}{ccc}0&0&0\\0&1&0\\0&0&-1\end{array}\right),
\eea

\bea
e_{1}=\left(\begin{array}{ccc}0&1&0\\0&0&0\\0&0&0\end{array}\right),\
e_{2}=\left(\begin{array}{ccc}0&0&0\\0&0&1\\0&0&0\end{array}\right),\
e_{{12}}=\left(\begin{array}{ccc}0&0&1\\0&0&0\\0&0&0\end{array}\right),
\eea
\bea
f_{1}=\left(\begin{array}{ccc}0&0&0\\1&0&0\\0&0&0\end{array}\right),\
f_{2}=\left(\begin{array}{ccc}0&0&0\\0&0&0\\0&1&0\end{array}\right),\
f_{{12}}=\left(\begin{array}{ccc}0&0&0\\0&0&0\\1&0&0\end{array}\right).
\eea
Quadratic Casimir operator:
\bea
C_2=\sum_{i=1}^3
(e_{i}f_{i}+f_{i}e_{i})+{2\over 3}
\left(\sum_{i=1}^2h_{i}^2+h_{1}h_{2}
\right)=\\=
2\sum_{i=1}^3f_{i}e_{i}+2\sum_{i=1}^2 h_{i}+
{2\over 3}\left(\sum_{i=1}^2h_{i}^2+h_{1}h_{2}\right).
\eea

{Highest weight representation:}

We consider the principal (spherical) series of
representations, induced by the
one-dimensional representations of the Borel subgroup.
Let us describe Borel-Weil realization.
The space of representation $V_\lb$ $\lb=(\lb_1,\lb_1)$
is functions of three real variable $x_1,x_2,x_{12}$,
which are  matrix elements of the $3*3$ upper triangular matrices, i.e.
coordinates on the biggest Schubert cell of the $SL(3)$ flag variety $G/B$.
Matrix elements are defined by integrals with the flat measure.
The action of the algebra is given by the differential operators:
%\label{bwb-sl3}
\bea
e_{1} = \frac{\partial}{\partial x_1}, \ \
e_{2} = \frac{\partial}{\partial x_2} +
         x_1 \frac{\partial}{\partial x_{12}},\ \
e_{{12}} = \frac{\partial}{\partial x_{12}},
\eea
$$
h_{1} =-2x_1{\d x_1}+x_2{\d x_2}-x_{12}{\d x_{12}}+\lb_1,
$$
$$
h_{2} = x_1{\d x_1}-2x_2{\d x_2}-x_{12}{\d x_{12}}+\lb_2,
$$
\bea
f_{1} = \lb_1x_1- x_1^2\frac{\partial}{\partial x_1} -
      x_1x_{12}\frac{\partial}{\partial x_{12}} -
      (x_{12} - x_1x_2)\frac{\partial}{\partial x_2},\\
f_{2} = \lb_2x_2+
         x_{12}\frac{\partial}{\partial x_1} -
          x_2^2\frac{\partial}{\partial x_2} ,  \nn\\
f_{{12}} = \lb_1x_{12}-\lb_2(x_1x_2-x_{12})-x_1x_{12}
{\d x_1}+x_2(x_1x_2-x_{12}){\d x_2}-x_{12}^2{\d x_{12}}.
\eea

Obviously that $Constants$ are vacuum (i.e. the highest weight)
vectors in such representation, i.e.
$e_i (Constant) = 0,~~~~h_i (Constant)=\lb_i (Constant)$.

\subsection{Whittaker function.}
%%%%%%%%%%%%%%%%%%%%%%%%%%%%%%%%%%%%%%%%%%%%%%%%%%%%%%%%%%%%%%%%%%%
%%%%%%%%%%%%%%%%%%%%%%%%%%%%%%%%%%%%%%%%%%%%%%%%%%%%%%%%%%%%%%%%%%%

{\Large Definition 3.1:}
The vector $|w>_{\lb}^{\mu_1,\mu_2} $
is called {\it Whittaker vector} in representation $V_\lb$,
if it is eigenvector for generators $e_1,e_2$ i.e.
$e_1 |w>_{\lb}^{\mu_1,\mu_2}=\mu_1 |w>_{\lb}^{\mu_1,\mu_2}$ and
$e_2 |w>_{\lb}^{\mu_1,\mu_2}=\mu_2 |w>_{\lb}^{\mu_1,\mu_2}$
for some constants $\mu_i$. Obviously $e_3 |w>_{\lb}^{\mu_1,\mu_2}=0$.
We will sometimes omit indexes $\lb,\mu$, if they remains unchanged
during the calculations.

{\lL 3.2:}
For given $\mu_i$ and the irreducible representation $V_\lb$~ ($\lb_i \ne 1,2,3,...$)
there exist a unique (up to scalar)
Whittaker vector $\w \in V_\lb$, which can be expressed as follows:
\bea
\w = exp (\mu_1 x_1+\mu_2 x_2)
%= \sum_{i=0} \frac{\mu^n f^n}{n! (\lb,n)} |vac>
\label{wh_sl3}
\eea
%where we denote $ (\lb,n)=\lb (\lb-1) ... (\lb-n+1)$.
The proof is obvious.

%\lL 3.3: in terms of generators $f_i$ Whittaker vector can
%be given by
%\bea
%|w>_\lb=\sum_{n=0}^{\infty} \sum_{i=0}^{n}
%\eea

\Rem in case of $sl(2)$ we also gave an expression for the
Whittaker vector $|w>$ in terms of the generator $f$, i.e.
$|w>= \sum_{i=0} \frac{\mu^n f^n}{n! (\lb,n)} |vac>$
one can see that such expression in case of $sl(3)$ is rather
complicated and unusable.

{\lL 3.3: } The dual Whittaker vector: $<w|$ is given by the formula:
\bea
<w|=(x_{12}-x_1x_2)^{-(\lb+2)}x_{12}^{-(\lb+2)}
e^{-\mu_2^L{x_1\over x_{12}}-\mu_1^L{x_2\over x_1x_2-x_{12}}}.
\eea
see \cite{Liouvil}.

{\lD 3.4:} The Whittaker function $W_{\lb}^{\mu_L,\mu_R}(\phi_1,\phi_2)$
is the function given by:
\bea
W_\lb^{\mu_L,\mu_R} (\phi_1,\phi_2)= <w|exp( \phi_1 h_1+\phi_2 h_2)|w>
\label{wf_sl3}
\eea
We will sometimes omit indexes $\lb,\mu$.

{\lPn 3.5:}
The Whittaker function satisfies the equation:
\bea\label{HI2-sl3}
&&
\left(\frac{2}{3}\left({\partial^2\over\partial\phi_1^2}+
{\partial^2\over\partial\phi_2^2}+{\partial^2\over\partial\phi_1\partial\phi_2}
\right)+{\d\phi_1}+{\d\phi_2}-
\mu_1^L\mu_1^Re^{\phi_2-2\phi_1}-\mu_2^L\mu_2^Re^{\phi_1-2\phi_2}
\right) W_{\lb}(\phi_1,\phi_2)=\nn\\
&&
(2\lb_1+2\lb_2+\frac{2}{3}(\lb_1^2+\lb_2^2+\lb_1\lb_2) W_{\lb}(\phi_1,\phi_2).
\eea
\Rem the proof is based on idea that this hamiltonian is the second order
Casimir operator for the $SL(3)$, and completely analogous to the
$SL(2)$ case. The function $W_{\lb}(\phi_1,\phi_2)$ is also eigenfunction
for the third order hamiltonian, which is the image of the third order Casimir
for the  $SL(3)$. So we see that  Whittaker function is the wave function
for the $SL(3)$ quantum Toda chain. These facts are standard so we are very
brief.

{\lPn 3.6 (Integral representation of the Whittaker function) (\cite{Liouvil})
:}
\bea \label{LWFG-sl3} W_\lb(\phi_1,\phi_2)=
e^{(\lb_1+1)\phi_1+(\lb_2+1)\phi_2}\int \ dx_1dx_2dx_{12}\
(x_{12}-x_1x_2)^{-\lb_2-2)}x_{12}^{-\lb_1-2)}\times\\\times
e^{-\mu_2^L{x_1\over x_{12}}-\mu_1^L{x_2\over x_1x_2-x_{12}}-
\mu_1^Rx_1e^{\phi_2-2\phi_1}-\mu_2^Rx_2e^{\phi_1-2\phi_2}}.
\eea

The theorem is simple corollary of the facts  that the invariant pairing
$<v|u>$ is realized as integral $\int$ with the flat measure, and
dual Whittaker vector $<w|$ is realized as
$(x_{12}-x_1x_2)^{(-\lb_2-2)}x_{12}^{(-\lb_1-2)}
exp(-\mu_2^L{x_1\over x_{12}}-\mu_1^L{x_2\over x_1x_2-x_{12}})$.

\subsection{Intertwining operators}
%%%%%%%%%%%%%%%%%%%%%%%%%%%%%%%%%%%%%%%%%%%%%%%%%%%%%%%%%%%%%%%%%%%
%%%%%%%%%%%%%%%%%%%%%%%%%%%%%%%%%%%%%%%%%%%%%%%%%%%%%%%%%%%%%%%%%%%

Denote by $V_{(1,0)}$ %($V_{(0,1)}$)
the first %(second)
fundamental representation of the $SL(3,\RR)$,
denote  by $|0>$ the highest weight vector, by $|1>$ the vector
$f_1|0>$, by $|2>$  the vector $f_2 f_1|0>=f_3|0>$, i.e.:
\bea
|0>=\left(\begin{array}{cc}1\\0\\0\end{array}\right),\
|1>=\left(\begin{array}{cc}0\\1\\0\end{array}\right),\
|2>=\left(\begin{array}{cc}0\\0\\1\end{array}\right),\
%$$ %
\eea
hence $e_1 |0>=0, e_2|0>=0, h_1 |0>=|0>, h_2 |0>=0$.

%Let us construct the  intertwining operators $V_{\lb+1} \to \vL\otimes V_1$

{\lPn 3.7:}
Let us denote by
$\Phi_{\lb} = \Phi_{\lb,(+1,0)} \oplus \Phi_{\lb,(-1,1)}
\oplus \Phi_{\lb,(0,-1)}$ the ismorphism $
V_{\lb+(1,0)} \oplus V_{\lb+(-1,1)} \oplus  V_{\lb+(0,-1)} \to \vL\otimes V_1$.
The operators $\Phi_{\lb,(+1,0)}$ and ${\Phi_{\lb,(+1,0)}}^{-1}$ in Borel-Weil
realization are given by the formulas:
\bea
&&\Phi_{\lb,(+1,0)}:x_1^{n_1}x_2^{n_2}x_3^{n_3}|0>_{\lb}\otimes |0> \to
x_1^{n_1}x_2^{n_2}x_3^{n_3}|0>_{\lb+(1,0)},\\
&&\Phi_{\lb,(+1,0)}:
x_1^{n_1}x_2^{n_2}x_3^{n_3}|0>_{\lb}\otimes |1> \to
x_1^{n_1+1}x_2^{n_2}x_3^{n_3}|0>_{\lb+(1,0)},\\
&&\Phi_{\lb,(+1,0)}:
x_1^{n_1}x_2^{n_2}x_3^{n_3}|0>_{\lb}\otimes |2> \to
x_1^{n_1}x_2^{n_2}x_3^{n_3+1}|0>_{\lb+(1,0)}.
\label{inter-sl3-bwb-easy}
\eea

\bea &&{(\Phi_{\lb,(+1,0)})}^{-1}:\nn\\
&&x_1^{n_1}x_2^{n_2}x_3^{n_3}|0>_{\lb+(1,0)} \to
\frac{1}{(\lb_1+1)(\lb_1 +\lb_2+2)} (
(n_1 n_2 x_3 +((\lb_1+2)n_3+n_2n_3)x_1 x_2)
x_1^{n_1-1}x_2^{n_2-1}x_3^{n_3-1}|0>_{\lb} \otimes
|2>+\nn\\
&&
(n_1(\lb_1+\lb_2+2-n_2)x_3
+n_3(\lb_2-n_2)n_3 x_2) x_1^{n_1-1}x_2^{n_2}x_3^{n_3-1})|0>_{\lb} \otimes |1>+\nn\\
&&
  (-n_1n_2 x_3^{2} +
n_3 (-\lb_2+n_2) x_1^{2}x_2^{2}+\nn \\
&&
+( (\lb_1+1)(\lb_1+\lb_2+2)-n_1(\lb_1+\lb_2+2)+n_1n_2
- n_3(\lb_1+2)-n_2n_3)x_1x_2x_3)
x_1^{n_1-1}x_2{n_2-1}x_3^{n_3-1}|0>_{\lb} \otimes |0>).
\label{intert-sl3-bwb}
\eea

This proposition is auxillary for us. The main formulas we need to obtain
the relations between the Whittaker wave functions are the formulas
from the next proposition. But to obtain them we need this statement.
Actually we use it rather small, only to fix the constants (see below).

To prove the formula \ref{inter-sl3-bwb-easy} one should only note
that the representation $V_{(1,0)}$ also have Borel-Weil realization and
the above intertwiner is given by the multiplication of functions from
the two representations. The proof of the formula  \ref{intert-sl3-bwb}
is more complicated, but it is standard representation theory reasoning,
and should be available in the literature, but we were not able to find it.
So we sketch the proof. Let us denote the three components of the
operator ${(\Phi_{\lb,(+1,0)})}^{-1}$ by $\Psi_0,\Psi_1, \Psi_2$, i.e.
${(\Phi_{\lb,(+1,0)})}^{-1} |v>= \Psi_0 (|v>)\otimes |0>
+\Psi_1(|v>)\otimes |1>+ \Psi_2(|v>)\otimes |2>$,
first simple, but useful step in proving   \ref{intert-sl3-bwb}
is that all the three operators $\Psi_i$
can be easily expressed one throw another,
so it's enough to find only one of them. Actually it is easy to see that, for
example:
$\Psi_1=[\Psi_2,f_2]$, $\Psi_0=[\Psi_2,f_3]$. Looking at the formula
\ref{intert-sl3-bwb} it is clear that operator $\Psi_2$ is much simpler than
the other two, so it is rather useful that we can express them through
the operator $\Psi_2$. Second step is to find $\Psi_2$. It goes as follows,
first let us note that: $\Psi_2 (x_1^{k_1}x_2^{k_2}x_3^{k_3} )=
\alpha_{k_1,k_2,k_3} x_1^{k_1-1}x_2^{k_2-1}x_3^{k_3} +
\beta_{k_1,k_2,k_3} x_1^{k_1}x_2^{k_2}x_3^{k_3-1}$. That easily follows
from the equalities:
$\Psi_2 h_1=h_1\Psi_2, \Psi_2 h_2=h_1 \Psi_2$ and
$e_1^{k_1+1} \Psi_2 (x_1^{k_1}x_2^{k_2}x_3^{k_3} )=0,
e_3^{k_1+1} \Psi_2 (x_1^{k_1}x_2^{k_2}x_3^{k_3} )=0$.
the second point is to use the commutation relation of  $e_2$ and $\Psi_2$
and to find  recurrence relations for
$\alpha_{k_1,k_2,k_3}, \beta_{k_1,k_2,k_3}$, they can be solved directly
after some work. Hence we find $\Psi_2$ and as explained above it follows
that we find $\Phi_1,\Phi_0$.

\Rem Some other useful formulas for the intertwiners we put to the Appendix.
Let us note that it is rather easy to find the formulas for the intertwiners
in Verma realization (see Appendix), not Borel-Weil, but the problem is that
expression of  the vector $x_1^{k_1}x_2^{k_2}x_3^{k_3}$ through the vectors
$f_1^{n_1}f_2^{n_2}f_3^{n_3}$ is rather complicated  and
untreatable.

\lPn 3.8:
The action of the intertwiner ${(\Phi_{\lb,(+1,0)})}$
on the Whittaker vector can be expressed as follows:
\bea
\label{int-1-sl3-whit}
&&
{(\Phi_{\lb,(+1,0)})}  |w>_{\lb}^{\mu_1, \mu_2}\otimes |0>
\to |w>_{\lb+(1,0)}^{\mu_1, \mu_2}=
exp(\mu_1x_1+ \mu_2 x_2)|0>_{\lb+(1,0)},\\
&&
{(\Phi_{\lb,(+1,0)})}  |w>_{\lb}^{\mu_1, \mu_2}\otimes |1>
\to %\frac{1}{\mu_1}
\bar{\Lambda}_{1} |w>_{\lb+(1,0)}^{\mu_1, \mu_2}=
x_1 exp(\mu_1x_1+ \mu_2 x_2)|0>_{\lb+(1,0)},\\
&&
{(\Phi_{\lb,(+1,0)})}  |w>_{\lb}^{\mu_1, \mu_2}\otimes |2>
\to
\frac{1}{\mu_1}(f_2+\bar{\Lambda}_{2}\mu_2 -\lb_2\bar{\Lambda}_{2}
|w>_{\lb+(1,0)}^{\mu_1, \mu_2}=x_{12} exp(\mu_1x_1+ \mu_2 x_2)|0>_{\lb+(1,0)}.
\eea
where
$ \bar{\Lambda}_{1}=\frac{2(\lb_1+1)+\lb_2 -2 h_1 -h_2}{3\mu_1}$ and $
\bar{\Lambda}_{2} = \frac{(\lb_1+1)+2\lb_2 - h_1 -2 h_2}{3\mu_2}
$

The action of the intertwiner ${(\Phi_{\lb,(+1,0)})}^{-1}$
on the Whittaker vector can be expressed as follows:
\bea
\label{int-2-sl3-whit}
{(\Phi_{\lb,(+1,0)})}^{-1} |w>_{\lb+1}^{\mu_1, \mu_2} \to
\frac{1}{(\lb_1+1)(\lb_1 +\lb_2+2)} (
\mu_1\mu_2 |w>_{\lb}^{\mu_1, \mu_2} \otimes |2>+\nn \\
\frac{1}{\mu_2}(\tilde {\Lambda}_2+2)|w>_{\lb}^{\mu_1, \mu_2} \otimes |1>+
\frac{1}{\mu_2}(f_1+ \frac{1}{\mu_2} (\tilde{\Lambda}_1^2-
(\lb_2-3) \tilde{\Lambda}_1)) |w>_{\lb}^{\mu_1, \mu_2} \otimes |0>).
\eea
where $\tilde {\Lambda}_1= \frac{2h_1+h_2 +\lb_1+2\lb_2}{3}$
and  $\tilde {\Lambda}_2= \frac{2h_2+h_1 +\lb_2+2\lb_1}{3}$

This proposition is crucial in this section. The relations for the
Whittaker wave functions (\ref{rais_1_sl3},\ref{bil_sl3}) easily follows from it.
Let us sketch the proof and explain the difference with $SL(2)$ case,
which is much more simple. We discuss the analogous reasonings
with more details in the next section in the case of $SL(n)$.
Let us  denote by $P_0,P_1,P_2$ the operators such that:
${(Phi_{\lb,(+1,0)})}^{-1} |w>_{\lb+1}=
\Psi_0 |w>_{\lb}\otimes |0>+\Psi_1 |w>_{\lb}\otimes |1>+
\Psi_2 |w>_{\lb}\otimes |2>.$ We will look for them as some polynomials
$P_i(h_i,f_i)$ from the generators $e_i$.
From the commutation relations with the operators $e_i$ and uniqueness
of the Whittaker vector it's easy
to see that the operator $P_2$ is some constant depending only on
$\mu_i$. To fix this constant we look on the formulas \ref{intert-sl3-bwb},
\ref{wh_sl3}
and see that $P_2=\mu_1\mu_2$. Also from the
 commutation relations with operators $e_i$ follows the following relations:
 $[P_1,e_1]=0$, $[P_1,e_2]=e_2+Constant$ so $A_1$ equals
$\frac{1}{\mu_2}\Lambda_2+C_1$, where $\Lambda_i$ are fundamental coroots,
i.e. elements of Cartan subalgebra such that $[\Lambda_i,e_j]=\delta_{i,j}$,
the  constant $C_1$ can also be fixed  by the formulas
\ref{intert-sl3-bwb},\ref{wh_sl3}. Analogously $A_0$ satisfy
the identities: $[P_0,e_2]=0$, $[P_0,e_1]=\frac{\alpha}{\mu_1}A_1e_1+
\beta P_1$, where $\alpha+\beta=1$.
the solution of this equation is given by the ansatz $P_0=\frac{1}{\mu_2}
(f_1+\frac{1}{\mu_1}(\Lambda_1^2+C_0\Lambda_1))$ the
constants $c_0, \alpha, \beta$ also can be fixed by the formulas
\ref{intert-sl3-bwb},\ref{wh_sl3}. So the necessary formulas
for $P_i$ are obtained.

Let us one more time emphasize the difference with the $SL(2)$ case,
where we had the formula for the Whittaker vector
$\w=\sum_{i=0} \frac{\mu^n f^n}{n! (\lb,n)} |vac>$. Hence the proposition
2.8, which is analogous to the proposition 3.8 above immediately follows
from the formula for the action of the intertwiner on the vector $f^n|vac>$.
In the present case the formula for $\w$ in terms of the generators
$f_i$ is too complicated to extract something from it, so the proof
of the proposition above was based mostly on the property of the
Whittaker vector: $e_i|w>=\mu_i|w>$.

\subsection{Relations for the Whittaker wave functions.}
%%%%%%%%%%%%%%%%%%%%%%%%%%%%%%%%%%%%%%%%%%%%%%%%%%%%%%%%%%%%%%%%%%%
%%%%%%%%%%%%%%%%%%%%%%%%%%%%%%%%%%%%%%%%%%%%%%%%%%%%%%%%%%%%%%%%%%%

%$$ $$
{$~$ }

{\lT (Raising operators) 3.9:} The following relations holds:
\bea
\label{rais_1_sl3}
&&
  W_{\lb+1}(\phi_1,\phi_2)=
 \frac{exp(\phi_1)}{\mu_1^{R}(\lb_1+1)(\lb_1+\lb_2+2)}
( \mu_1^R exp(-2\phi_1 +\phi_2) +
\frac{1}{\mu_2^R}(
\left(\frac{2\partial_{\phi_1}+\partial_{\phi_2}+\lb_1+2\lb_2}{3}\right)^2\nn \\
&&
-(\lb_2-3)\frac{2\partial_{\phi_1}+\partial_{\phi_2}+\lb_1+2\lb_2}{3}
-2(\lb_2-1))) W_{\lb}(\phi_1,\phi_2).
\eea

 This proposition follows from the \ref{int-1-sl3-whit},\ref{int-2-sl3-whit}
in a way completely analogous to the $SL(2)$ case  (proposition 2.11),
so we omit the details.

\lT (Bilinear relations) 3.10:
\bea
\label{bil_sl3}
&&
(
(\frac{1}{\mu_2^{\prime L }} ( \mu_2^{\prime R } exp(-\phi_2 +2\phi_1) +
(
\left(\frac{-\partial_{\phi_2^{\prime}}-2\partial_{\phi_1^{\prime}}+
(\nu_2+1)+ 2\nu_1}{3\mu_1^{\prime L}}\right)^2 \nn \\
&&
-(\nu_1)\frac{-\partial_{\phi_2^{\prime}}-2\partial_{\phi_1^{\prime}}+
(\nu_2+1)+ 2\nu_1}{3\mu_1^{\prime L}})
-
\frac{-2\partial_{\phi_1}-\partial_{\phi_2}+2(\lb_1+1)+ \lb_2}{3\mu_1^{L}}
\frac{-2\partial_{\phi_2^{\prime}}-\partial_{\phi_1^{\prime}}+
2(\nu_2+1)+ \nu_1}{3\mu_2^{\prime L}}\nn\\
&&
+
(\frac{1}{\mu_1^{L}} ( \mu_1^{ R} exp(-\phi_1 +2\phi_2) +
(
\left(\frac{-\partial_{\phi_1^{ }}-2\partial_{\phi_2^{ }}+
(\nu_1+1)+ 2\nu_2}{3\mu_2^{  L}}\right)^2 \nn\\
%\eea
%\bea
&&
-(\nu_2)\frac{-\partial_{\phi_1^{ }}-2\partial_{\phi_2^{ }}+
(\nu_1+1)+ 2\nu_2}{3\mu_2^{  L}})
)
W_{\lb+1}^{\mu_1^{L},\mu_2^{L},\mu_1^{R},\mu_2^{R}}(\phi_1,\phi_2)
W_{\nu+1}^{\mu_1^{\prime L},\mu_2^{\prime L},\mu_1^{\prime R},\mu_2^{\prime R}}
(\phi_1^{\prime},\phi_2^{\prime})
=\nn\\
%%%%%%%%%%%%%%%%%%%%%%%%%%%%%%%%%%%%%%%%%%%%%%%%%%%%%%%%%%%%%%%%%%%%%%%%
%%%%%%%%%%%%%%%%%%%%%%%%%%%%%%%%%%%%%%%%%%%%%%%%%%%%%%%%%%%%%%%%%%%%%%%%
%%%%%%%%%%%%%%%%%%%%%%%%%%%%%%%%%%%%%%%%%%%%%%%%%%%%%%%%%%%%%%%%%%%%%%%%
%%%%%%%%%%%%%%%%%%%%%%%%%%%%%%%%%%%%%%%%%%%%%%%%%%%%%%%%%%%%%%%%%%%%%%%%
&&
 \frac{1}{(\lb_1+1)(\lb_1+\lb_2+2)}\frac{1}{(\nu_2+1)(\nu_1+\nu_2+2)}
%\eea
%\bea
(
(\frac{1}{\mu_1^{R}} ( \mu_1^L exp(-2\phi_1 +\phi_2) +
\frac{1}{\mu_2^R}(
\left(\frac{2\partial_{\phi_1}+\partial_{\phi_2}+\lb_1+2\lb_2}{3}\right)^2 \nn\\
&&
-(\lb_2-3)\frac{2\partial_{\phi_1}+\partial_{\phi_2}+\lb_1+2\lb_2}{3}
-2(\lb_2-1)))
)\mu_1^{\prime R} \mu_2^{\prime R} %\nn\\
%&&
+
\frac{(\partial_{\phi_1}+2\partial_{\phi_2}+2\lb_1+\lb_2)}{3 \mu_2}
\frac{(\partial_{\phi_1^{\prime}}+2\partial_{\phi_1}+2\nu_1+\nu_1)}{3 \mu_1}
-
\nn\\
%\eea
%\bea
&&
-\mu_1^{} \mu_2^{}
(\frac{1}{\mu_2^{\prime R}} ( \mu_2^{\prime L} exp(-2\phi_2^{\prime} +\phi_1^{\prime}) +
\frac{1}{\mu_1^{\prime R}}(
\left(\frac{2\partial_{\phi_2^{\prime}}+\partial_{\phi_1^{\prime}}+\nu_2+2\nu_1}{3}\right)^2 \nn\\
%\eea
%\bea
&&
-(\nu_1-3)\frac{2\partial_{\phi_2^{\prime}}+\partial_{\phi_1^{\prime}}+\nu_2+2\nu_1}{3}
-2(\nu_1-1)))
)
W_{\lb}^{\mu_1^{L},\mu_2^{L},\mu_1^{R},\mu_2^{R}}(\phi_1,\phi_2)
W_{\nu}^{\mu_1^{\prime L},\mu_2^{\prime L},\mu_1^{\prime R},\mu_2^{\prime R}}
(\phi_1^{\prime},\phi_2^{\prime}).
\eea

 The prooof  is also completely analogous to the $SL(2)$ case, so
we refer the reader to the previous section.

One can obtain the nonlinear equations from this bilinear one,
as it was done in the $SL(2)$ case, but we omit these considerations
due to their length.

%%%%%%%%%%%%%%%%%%%%%%%%%%%%%%%%%%%%%%%%%%%%%%%%%%%%%%%%%%%%%%%%%%%
%%%%%%%%%%%%%%%%%%%%%%%%%%%%%%%%%%%%%%%%%%%%%%%%%%%%%%%%%%%%%%%%%%%
%%%%%%%%%%%%%%%%%%%%%%%%%%%%%%%%%%%%%%%%%%%%%%%%%%%%%%%%%%%%%%%%%%%
%%%%%%%%%%%%%%%%%%%%%%%%%%%%%%%%%%%%%%%%%%%%%%%%%%%%%%%%%%%%%%%%%%%
%%%%%%%%%%%%%%%%%%%%%%%%%%%%%%%%%%%%%%%%%%%%%%%%%%%%%%%%%%%%%%%%%%%
%%%%%%%%%%%%%%%%%%%%%%%%%%%%%%%%%%%%%%%%%%%%%%%%%%%%%%%%%%%%%%%%%%%
%%%%%%%%%%%%%%%%%%%%%%%%%%%%%%%%%%%%%%%%%%%%%%%%%%%%%%%%%%%%%%%%%%%
%%%%%%%%%%%%%%%%%%%%%%%%%%%%%%%%%%%%%%%%%%%%%%%%%%%%%%%%%%%%%%%%%%%
\section{$ SL(n)$ Toda Chain.} \label{sl(n)}
%%%%%%%%%%%%%%%%%%%%%%%%%%%%%%%%%%%%%%%%%%%%%%%%%%%%%%%%%%%%%%%%%%%
%%%%%%%%%%%%%%%%%%%%%%%%%%%%%%%%%%%%%%%%%%%%%%%%%%%%%%%%%%%%%%%%%%%
%%%%%%%%%%%%%%%%%%%%%%%%%%%%%%%%%%%%%%%%%%%%%%%%%%%%%%%%%%%%%%%%%%%
%%%%%%%%%%%%%%%%%%%%%%%%%%%%%%%%%%%%%%%%%%%%%%%%%%%%%%%%%%%%%%%%%%%
%%%%%%%%%%%%%%%%%%%%%%%%%%%%%%%%%%%%%%%%%%%%%%%%%%%%%%%%%%%%%%%%%%%
%%%%%%%%%%%%%%%%%%%%%%%%%%%%%%%%%%%%%%%%%%%%%%%%%%%%%%%%%%%%%%%%%%%

\subsection{Notations.}
Algebra $sl(n)$ is completely given by the generators $e_i,h_i,f_i$, $i=1,...,n-1$
and commutation relations:
\bea\label{commrel1}
%\phantom{fhg}
[e_{ i},h_{j}]= -A_{ij}e_{\pm,i},\ \
[f_{ i},h_{j}]= A_{ij}e_{\pm,i},\ \
[e_{i},f_{j}]=
\delta_{ij}h_{j},\ \ i,j=1,\ldots, N-1,
\eea
and the Serre relations
\bea\label{Serre}
\phantom{fhg}\hbox{ad}_{T_{\pm i}}^{1-A_{ij}}\left(T_{\pm j}\right)=0,
\eea
where $\hbox{ad}_x^k(y)\equiv \underbrace{[x,[x,...,[x,y]..]]}_{k\ \hbox{
times}}$.

Where $A_{ij}$ is Cartan matrix, which is for the algebra $SL(N)$
is equal to
$$
\pmatrix
          2&-1&0 &0 &...&0\cr
                       -1&2 &-1&0 &...&0\cr
                        0&-1&2 &-1&...&0\cr
              \vdots& \vdots& \ddots &\ddots& \ddots &\vdots\cr
                        0& 0& 0&0 & -1&2
\endpmatrix
$$

Quadratic Casimir operator is
\bea
C_2=\sum_{\alpha\in\Delta}
e_{\alpha}f_{\alpha}+\sum_{ij}^{N-1}A_{ij}^{-1}h_{i}h_{j},
\eea
where the first sum goes over all (positive and negative) roots.

In the case of generic $SL(N)$ group, one can define the
(right) regular representation only in general terms of the group acting on the
space of the algebra of functions:
\bea\label{rrrep}
\pi_{reg}(h) f(g)=f(gh).
\eea
Therefore, we use from now on mostly group (not algebra) terms. Still, we can
restrict the space of functions to the irreducible representations in the
generic situation. For doing this, we consider the representation denoted
by $V_{\lb}$ induced by
one-dimensional representations of the Borel subgroup. That is, we reduce the
space of all functions to the functions satisfying the following covariance
property:
\bea\label{irrep}
f_{\lambda}(bg)=\chi _{\lambda}(b)f_{\lambda}(g),
\eea
where $b$ is an element of
the Borel subgroup of lower-triangle matrices and $\chi _{\lambda}$ is the
character of the Borel subgroup of the form:
\bea \label{rep}
\chi_{\lambda}(b)=\prod _{i=1}^{N-1} \mid
b_{ii} \mid ^{ (\lambda-\rho)  {\mbox{\bf $e$ }}_i}
(\mbox{sign}  b_{ii} )^{\epsilon _i},
\eea
where $\epsilon_i$ are equal to either 0 or 1. For the sake of simplicity, we
consider the representations with all these sign factors to be zero although
other cases can be also easily treated. The representation constructed belongs
to the principal (spherical) series.

Thus, our representation is given by
restricting the space of functions to the functions defined on the coset
$B\backslash G$, the biggest cell of  which, in turn, may be identified with the strictly
upper-triangular matrices $N_+$.
We will denote by $x_{i,j}$ the matrix elements of $n*n$ matrices, and
discussion above gives us realization of the representation of the $SL(n)$
in the space of functions of $x_{i,j}$ , $1\leq i<j\leq n$.
The explicit formulas can be found for example in \cite{Etingof, Frenkel-Feigin}.

\subsection{Whittaker function.}
%%%%%%%%%%%%%%%%%%%%%%%%%%%%%%%%%%%%%%%%%%%%%%%%%%%%%%%%%%%%%%%%%%%
%%%%%%%%%%%%%%%%%%%%%%%%%%%%%%%%%%%%%%%%%%%%%%%%%%%%%%%%%%%%%%%%%%%

{\lD 4.1:}
The vector $|w>_{\lb}^{\mu_1,...,\mu_{n-1}} $
is called {\it Whittaker vector} in representation $V_\lb$,
if it is eigenvector for generators $e_1,...,e_{n-1}$ i.e.
$e_i |w>_{\lb}=\mu_i |w>_{\lb}$
for some constants $\mu_i$. Obviously $[e_i,e_j] |w>_{\lb}^{\mu_1,\mu_2}=0$.
We will sometimes omit indexes $\lb,\mu$, if they remains unchanged
during the calculations.

{\lL 4.2:}
For given $\mu_i$ and irreducible representation $V_\lb$~
%($\vL$ is not finite-dimensional - dly nix toze proxodit )
there exist a unique (up to scalar)
Whittaker vector $\w \in V_\lb$, which can be expressed as follows:
\bea
\w = exp (\mu_1 x_{1,2}+\mu_2 x_{2,3}+...+\mu_{n-1} x_{n-1,n})
%= \sum_{i=0} \frac{\mu^n f^n}{n! (\lb,n)} |vac>
\label{wh_sln}
\eea

The proof is obvious due to explicit formulas for generators $e_i$
in Borel-Weil realization. $e_i$ acts on such vector
as $\partial_{x_{i,i+1}}$, see  \cite{Liouvil} for details.
%where we denote $ (\lb,n)=\lb (\lb-1) ... (\lb-n+1)$.

{\lD 4.3:} The Whittaker function $W_{\lb}^{\mu^L_i,\mu^R_j}
(\phi_1,\phi_2,...,\phi_{n-1})$
is the function given by:
\bea
W_\lb^{\mu^L_i,\mu^R_j} (\phi_1,\phi_2)=
<w|exp( \phi_1 h_1+\phi_2 h_2+...+\phi_{n-1}h_{n-1})|w>
\label{wf_sln}
\eea
We will sometimes omit indexes $\lb,\mu$.

{\lPn 4.4:}
The Whittaker function satisfies the equation:
\bea\label{HI2-sln}
\left(
\sum_{ij} A_{ij}^{-1}\partial_{\phi_i}\partial_{\phi_j}+
2\sum_{ij} A^{-1}_{ij} \partial_{\phi_j}
-
2\sum_i\mu^L_i\mu_i^R exp(\sum_j A_{i,j} \phi_j)
\right) W_{\lambda}^{\mu^L_i,\mu^R_j}(\phi_i)=
(\lambda^2-\rho^2)W_{\lambda}^{\mu^L_i,\mu^R_j}(\phi_i).
\eea
where $\rho$ is one half of sum of simple roots.

\Rem the proof is based on the idea, that this hamiltonian is the second order
Casimir operator for the $SL(n)$, and completely analogous to the
$SL(2)$ case. The function $W_{\lb}(\phi_i)$ is also eigenfunction
for the higher  order hamiltonians, which are the images of the higher
 order Casimirs
for the  $SL(n)$. Hence, we see that Toda is completely integrable system and
the  Whittaker function is the wave function
for all Hamiltonians of the $SL(n)$ quantum Toda chain. These facts are standard so we are very
brief (see \cite{Liouvil}).

{\lPn 4.5 (Integral representation of the Whittaker function)
(\cite{Liouvil}) :}
\bea
\label{LWFG-sln}
W_\lb(\phi_i)=
e^{-\sum_i \lambda_i\phi_i}
\int _{X=B\backslash G} \prod_{i<j}dx_{ij}\prod_{i=1}^{N-1}
\Delta_i^{-(\sum_k A_{k,i}\lambda_k+1)}(xS^{-1})\times
%\\ \times
e^{ \mu_{i}^Rx_{i,i+1}e^{(\sum_k A_{k,i}\lambda_k)} -\mu_{N-i}^L
\frac{\Delta _{i,i+1}(xS^{-1})}{\Delta _{i}(xS^{-1})}}.
\eea
where $\Delta_i$ - is i-th principal minor of matrix $x_{i,j}$,
$\Delta_{i,i+1}$ - determinant $(n-1)*(n-1)$ submatrix, which obtained
by interchanging $n-1$  and $n$ column. Matrix $S$ is given by
$S_{i,j}=\delta_{i+j,n}$.

The theorem is  corollary of the facts  that the invariant pairing
$<v|u>$ is realized as integral  with the flat measure, and
the realizations of the Whittaker and the dual Whittaker vectors.
The most nontrivial part is to find dual
Whittaker vector. We refer to \cite{Liouvil} for detailed exposition.

\subsection{Intertwining operators}
%%%%%%%%%%%%%%%%%%%%%%%%%%%%%%%%%%%%%%%%%%%%%%%%%%%%%%%%%%%%%%%%%%%
%%%%%%%%%%%%%%%%%%%%%%%%%%%%%%%%%%%%%%%%%%%%%%%%%%%%%%%%%%%%%%%%%%%

Denote by $V_{(1,0,...,0)}$ %($V_{(0,1)}$)
the first %(second)
fundamental representation of the $SL(n,\RR)$, which is just the standard
action of the $n*n$ matrices on ${\CC}^n$,
denote  by $|0>$ the highest weight vector, by $|1>$ the vector
$f_1|0>$, by $|2>$  the vector $f_2 f_1|0>=f_3|0>$, etc.

To construct relations for wave functions, one needs to express the
action of the intertwiner on Whittaker vector, in terms of the generators
$f_i,h_i$.
For example the crucial formula for us in the case
of $SL(3)$ was:

\bea
|w>_{\lb+1}^{\mu_1, mu_2} \to
\frac{1}{(\lb_1+1)(\lb_1 +\lb_2+2)} (
\mu_1\mu_2 |w>_{\lb}^{\mu_1, mu_2}
\otimes |2>+\\
\frac{1}{\mu_2}(\tilde {\Lambda}_2+2)|w>_{\lb}^{\mu_1, mu_2}
\otimes |1>+
\frac{1}{\mu_2}(f_1+ \frac{1}{\mu_2} (\tilde{\Lambda}_1^2-
(\lb_2-3) \tilde{\Lambda}_1)) |w>_{\lb}^{\mu_1, mu_2}
\otimes |0>).
\eea
Where $\tilde{\Lambda}_i$ were fundamental coweights minus some constants.

We need analogous formulas for the case of $SL(n)$.
We are unable to find them. We will present only the equations
solutions to which will give necessary formulas and first steps
towards  their solution.

Let us consider the most simple intertwiner
$\Phi: V_{\lb+(1,0,...,0)} \to V_\lb\otimes V_{(1,0,...,0)}$.
We need to find such polynomials $P_{j}(f_i,h_i)$, $0\leq j<n$ such that:
\bea
\Phi |w>_{\lb+(1,0,...,0)} = \sum_{j=0}^{n-1}
P_j(f_i,h_i) |w>_{\lb}\otimes |j>
\eea

{ \lPn 4.6}
The following relations holds :
\bea
e_k P_j(f_i,h_i) |w>_{\lb} + \delta_{j,k-1} P_k(f_i,h_i) |w>_{\lb}=
\mu_k P_j(f_i,h_i)|w>_{\lb}
\eea
The proof is obvious.

{\lC 4.7:}
One can look for $P_j(f_i,h_i) $ as an elements of the $U(sl(n))$,
satisfying the following relations:
\bea
\label{recurrent-sln}
[e_k, P_j(f_i,h_i)]=0, \mbox{ for } k\neq j+2\\ ~
[e_k, P_{k-2}(f_i,h_i)]=% ( \mbox{ modula } e_k=\mu_k) =
\alpha_k P_{k-1} (f_i,h_i)e_{k}+\beta_k P_{k-1} (f_i,h_i),
\mbox{ where } \alpha_k+\beta_k=1.
\eea

This allows one to organize itterative search of polynomials $P_j$:
if $P_j$ already found then $P_{j-1}$  is defined by conditions:
$
[e_k, P_{j-1}(f_i,h_i)]  = 0 $ for $k \neq j+1$ and
$ [e_{j+1}, P_{j-1}(f_i,h_i)]
=\alpha_k P_{k-1} (f_i,h_i)e_{k}+\beta_k P_{k-1} (f_i,h_i)$.

We can conclude from the considered examples that the highest degree
part of the $P_{j-1}$  can be defined uniquely from this relations.
The problem is that we do not know the way of solving such reccurence
relations.
Let us present the first 3 polynomials $P_j$ found by hands.

{\lPn 4.8:}
Polynomials $P_{n-1}, P_{n-2}, P_{n-3}$  are given by the formulas:
\bea
&& P_{n-1}= \mu_1 \mu_2 ... \mu_{n-1}, \ \
P_{n-2}= (\Lambda_{n-1} + Constant) \mu_1\mu_2 ... \mu_{n-2},\\
&&
P_{n-3}= (\mu_{n-2} f_{n-2}+ \mu_{n-3} f_{n-3}+ .... +\mu_{1} f_{1}+\nn\\
&&+
(\Lambda_{n-2}^2-\Lambda_{n-2}\Lambda_{n-3})+
(\Lambda_{n-3}^2-\Lambda_{n-3}\Lambda_{n-4})+...+(\Lambda_1)^2) \mu_1\mu_2 ... \mu_{n-3}
+ F_1({\Lambda_i}).
\label{first_pol}
\eea
where
$\Lambda_i $ - fundamental coweights, $F_1$ - some polynomial of degree 1.

The proof is straightforward check that given polynomials
satisfy relations (\ref{recurrent-sln}) and one can see that they are the only
with this property.

So we see polynomials $P_{n-1-j}$ are of degree j. It's easy to guess
that it also holds true for $j>3$.
It's easy to understand that
to write the simplest raising relation one should find $P_0$,
which should be some polynomial of degree $n-1$.
Our first hope to do that was to find some recurrence formula:
$P_{j-1}=A_j P_j$, but we did not succeed, looking at \ref{first_pol}
it seems that such formula may not exist.
Another idea is that may be $P_j$ are somehow related to the
Casimirs of subalgebras $SL(n-j)$, which one guess looking at
\ref{first_pol} and to commutation relation: $[e_k, P_j]=0$ for
$k\neq j+2$.

Let us also note that it is possible to write intertwiner
$V_{\lb+(1,0,...,0)} \to V_\lb\otimes V_{(1,0,...,0)}$
in Borel-Weil realization so it is possible to find the
expression for
 $ P_j $ as some polynomial of operators of multiplication on
 monomials $x_{i,j}$ but the problem
is then express such operators via the polynomial of the  $f_i,h_i$.

\section{Concluding remarks.} \label{conclud}

  Let us formulate the main points of our work. First,  we recalled
the representation theory approach to integrable systems, i.e.
we recalled  that appropriate matrix elements in irreducible representations
are the wave functions of integrable systems, Casimirs of the group
are  the Hamiltonians, the integral representation can be obtained
by realizing the representation of the group in space of functions,
where the pairing is given by the integral. Second,  we demonstrated in our
paper that different relations between the wave functions, can be obtained
from considerations of intertwining operators between  above representations.
We saw that raising, bilinear etc. relations can be obtained this way.

  In this paper we did not succeed in obtaining the general formula
for the case of $SL(n)$, it seems such formula can be quite complicated.
We only found some recurrence relations, whose solution leads to desired
formula. It seems that one should guess some ansatz for the desired
polynomial, something like determinant or trace of some matrix of the
generators of the algebra, like the formula for Casimirs.

  Group theory approach to integrable systems is very powerful,
it allows to give transparent explanations for different properties
of the system. So let us mention some questions about integrable systems
which should be interpreted by the group theory methods.
It will be interesting to understand the relation
between our approach to raising relations and approach based on Hecke
algebra (see \cite{Rais}). It is well-known that Bessel functions,
which are very close to  the $SL(2)$ Toda wave functions,
satisfy the following properties:

$ \sum_{n=-\infty}^{\infty} J_n(x) t^n=exp(\2x(t-\frac{1}{t})) $
   ~~~~~~~ and hence $J_{-n}(x)=(-1)^n J_{n}(x)$.

It will be interesting to find analogues formulas for the wave functions
of $SL(n)$ Toda chain and find their group theory interpretation.
In the studding integrable systems recently appeared such new important
ideas as shift operator and Dunkl operator \cite{Dunkl},  another new important
concept is Cherednik-Matsuo correspondence (see Cherednik in \cite{S}) between
solutions of integrable systems and solutions of Knizhnik-Zamolodchikov
equation it seems to be interesting
to understand their group theory meaning.

 Let us note the most straight-forward generalization of relations
obtained here is their transfer to the other groups:
instead of $SL(2,R)$ we plan to consider the  affine version
 $\hat{ SL}(2) $ related to periodic Toda chain,
quantum analogue $U_q( SL(2))$ related to Mcdonald polynomilas and
q-special functions, and may be finite field
analog $ SL(2,F_q)$, which may be related to the questions
considered in \cite{DL}.

\section{Acknowledgments.}

  The author is especially  grateful to S. Kharchev, for excellent
explanations and discussions, and sharing with the author his ideas,
which led to writing this paper. We are also grateful to A. Mironov
for the significant help during working on text, to A. Morozov
for initiating our work in this direction, and to our scientific advisor
S. Khoroshkin for constant support. The work was supported by
INTAS -93-10183 Ext., %(Khoroshkin)
Russian President's grant 96-15-96939
%INTAS-97-0103, %(Rosly)
and RFBR-98-0100344. %Khoroshkin
%????????? (Morozov).

$$ $$

{\Large Appendix.}

{$ ~ $}

In this Appendix we collect some more formulas on intertwining operators
for $SL(3)$. We use notations of section 3.

{\lPn:}
The representation  $ \vL\otimes V_1   $ is isomorphic to the :
$V_{\lb+(1,0)} \oplus V_{\lb+(-1,1)} \oplus  V_{\lb+(0,-1)} $.
The highest weight vectors in $ \vL\otimes V_1   $  are given by the
formulas:
\bea
|0>_{\lb}\otimes |0> - { \mbox{ weight }\lb+(1,0)} \\
f_1|0>_{\lb}\otimes |0>- \lb_1 |0>_{\lb}\otimes |1>- { \mbox{ weight }\lb+(-1,1)} \\
\lb_2(\lb_1+\lb_2+1)|0>_{\lb}\otimes |2>
-(\lb_1+\lb_2+1)f_2|0>_{\lb}\otimes |1> +
(f_1 f_2+\lb_2f_3) |0>_{\lb}\otimes |0> - { \mbox{ weight } \lb+(0,-1)}
\eea
The intertwiner
$\Phi_{\lb} = \Phi_{\lb,(+1,0)} \oplus \Phi_{\lb,(-1,1)}
\oplus \Phi_{\lb,(0,-1)} :
V_{\lb+(1,0)} \oplus V_{\lb+(-1,1)} \oplus  V_{\lb+(0,-1)} \to \vL\otimes V_1$
 in "Verma" basis is given by the formulas:
\bea
&&
\Phi_{\lb,(+1,0)} : f_1^{n_1}f_2^{n_2}f_3^{n_3}|0>_{\lb+(1,0)} \to
f_1^{n_1}f_2^{n_2}f_3^{n_3}|0>_{\lb}\otimes |0> +
n_1 f_1^{n_1-1}f_2^{n_2}f_3^{n_3}|0>_{\lb}\otimes |1> +
n_3 f_1^{n_1}f_2^{n_2}f_3^{n_3-1}|0>_{\lb}\otimes |2>, \\
&&
\Phi_{\lb,(-1,1)} : f_1^{n_1}f_2^{n_2}f_3^{n_3}|0>_{\lb+(-1,1)} \to
(f_1^{n_1+1}f_2^{n_2}f_3^{n_3}-n_2 f_1^{n_1}f_2^{n_2-1}f_3^{n_3+1})
|0>_{\lb}\otimes |0> + \nn \\
&&      ((-\lb_1+n_1)f_1^{n_1}f_2^{n_2}f_3^{n_3}-n_1n_2f_1^{n_1-1}f_2^{n_2-1}f_3^{n_3+1})
|0>_{\lb}\otimes |1> +
(n_2(\lb_2-n_3)f_1^{n_1}f_2^{n_2-1}f_3^{n_3}+
n_3f_1^{n_1+1}f_2^{n_2}f_3^{n_3-1})|0>_{\lb}\otimes |2>, \\
&&
\Phi_{\lb,(0,-1)} : f_1^{n_1}f_2^{n_2}f_3^{n_3}|0>_{\lb+(0,-1)} \to
(f_1^{n_1+1}f_2^{n_2+1}f_3^{n_3}+(\lb_2-n_2)
f_1^{n_1}f_2^{n_2}f_3^{n_3+1)}|0>_{\lb}\otimes |0> + \nn\\
&&     (n_1(\lb_2-n_2)f_1^{n_1-1}f_2^{n_2}f_3^{n_3+1}+
(n_1-(\lb_1+\lb_2+1))f_1^{n_1}f_2^{n_2+1}f_3^{n_3}) |0>_{\lb}\otimes |1>+
\nn\\
&&((\lb_2-n_2)(n_3-(\lb_1+\lb_2+1))f_1^{n_1}f_2^{n_2}f_3^{n_3)}+ n_3
f_1^{n_1+1}f_2^{n_2+1}f_3^{n_3-1}) |0>_{\lb}\otimes |2> .
\eea

The formula \ref{intert-sl3-bwb} for the intertwiner
${(\Phi_{\lb,(+1,0)})}$ in Borel-Weil realization
can be rewitten in more invariant terms as follows:
\bea
&&
{(\Phi_{\lb,(+1,0)})}^{-1} f(x_1, x_2,x_{12})|0>_{\lb+(1,0)} \to
\frac{1}{(\lb_1+1)(\lb_1 +\lb_2+2)}(
(\partial_{x_1} \partial_{x_2} +((\lb_1+2) + x_2 \partial_{x_2} )
\partial_{x_{12}} ) f(x_1, x_2,x_{12})|0>_{\lb} \otimes |2>+\nn \\
&&
( (\lb_1+\lb_2+2-x_2\partial_{x_2})\partial_{x_1}+
(\lb_2-x_2^2\partial_{x_2})\partial_{x_{12}})
f(x_1, x_2,x_{12})|0>_{\lb} \otimes |1>+\nn \\
&&
(-x_{12}\partial_{x_1})\partial_{x_{2}}+
x_1x_2(-\lb_2+x_2\partial_{x_{2}})\partial_{x_{12}}
+((\lb_1+1)(\lb_1+\lb_2+2) -(\lb_1+\lb_2+2)x_1\partial_{x_{1}}+\nn \\
&&
x_1x_2\partial_{x_{1}}\partial_{x_{2}}-(\lb_1+2)x_{1}\partial_{x_{1}}-
x_2x_{12}\partial_{x_{2}}\partial_{x_{12}}
)f(x_1, x_2,x_{12})|0>_{\lb} \otimes |0>).
\eea

{\lPn:}
In Borel-Weil realization the action of the intertwiner on the Whittaker
vector can be  written as follows:
\bea
&&
{(\Phi_{\lb,(+1,0)})}^{-1}: exp(\mu_1 x_1+ \mu_2 x_2)|0>_{\lb+1} \to
\frac{1}{(\lb_1+1)(\lb_1 +\lb_2+2)} (
\mu_1\mu_2  exp(\mu_1 x_1+ \mu_2 x_2)|0>_{\lb}\otimes |2>+\nn \\
&&
(\mu_1(\lb_1+\lb_2+2)-\mu_1\mu_2 x_2)
exp(\mu_1 x_1+ \mu_2 x_2)|0>_{\lb}\otimes |1>+\nn \\
&&
(\mu_1\mu_2(x_1x_2-x_{12})-\mu_1(\lb_1+\lb_2+2)x_1+
((\lb_1+1)(\lb_1+\lb_2+2))
exp(\mu_1 x_1+ \mu_2 x_2)|0>_{\lb}\otimes |0>).
\eea

\end{document}